\documentclass[12pt]{amsart}
\usepackage{amssymb}
\usepackage{amsmath}
\usepackage{epsfig}
\usepackage{graphics}

\font\tenmath=msbm10 \font\sevenmath=msbm7 \font\fivemath=msbm5
\newfam\mathfam \textfont\mathfam=\tenmath
\scriptfont\mathfam=\sevenmath \scriptscriptfont\mathfam=\fivemath

\def \\ { \cr }

\newcommand{\F}{{\mathcal F}}

\newcommand{\A}{{\mathcal A}}

\newcommand{\E}{{\mathcal E}}

\newcommand{\RR}{{\mathbb R}}
\newcommand{\CC}{{\mathbb C}}
\newcommand{\NN}{{\mathbb N}}
\newcommand{\ZZ}{\mathbb Z}

\newcommand{\EE}{{\mathbb E}}
\newcommand{\PP}{{\mathbb P}}

\newcommand{\vip}{\vskip0.15cm}

\newcommand{\indiq}{{\rm 1 \hskip-4pt 1}}

\def\P{{\mathcal P}}
\def\S{{\mathcal S}}

\newcommand{\xs}{X_{\sigma}}

\numberwithin{equation}{section}

\newtheorem{theo}{\indent Theorem}[section]
\newtheorem{prop}[theo]{\indent Proposition}
\newtheorem{con}[theo]{\indent Conjecture}
\newtheorem{rem}[theo]{\indent Remark}
\newtheorem{lem}[theo]{\indent Lemma}
\newtheorem{defin}[theo]{\indent Definition}
\newtheorem{cor}[theo]{\indent Corollary}

\theoremstyle{remark}
\newtheorem{remark}[theo]{Remark}

\newcommand{\sub}{{\sigma}} 
\newcommand{\dd}{d} 
\newcommand{\projinv}{\Pi_1} 

\begin{document}

\parindent = 0cm

\title[Deviation for substitution dynamical systems]{Deviation of ergodic averages for substitution dynamical systems with eigenvalues of modulus one}

\author{Xavier Bressaud}
\address{Universit\'e Paul Sabatier, Institut de Math\'ematiques de Toulouse, 118 route de Narbonne, F-31062 Toulouse Cedex 9, France}
\email{bressaud@math.univ-toulouse.fr}

\author{Alexander I. Bufetov}
\address{Rice University, Houston, TX, USA and the Steklov Institute of Mathematics, Moscow, Russian Federation} \email{aib1@rice.edu}

\author{Pascal Hubert}
\address{Laboratoire Analyse, Topologie et Probabilit\'es,
Case cour A,
Facult\'e des Sciences de  Saint-Jer\^ome,
Avenue Escadrille Normandie-Niemen,
13397 Marseille Cedex 20, France.} \email{hubert@cmi.univ-mrs.fr}

\subjclass{Primary: 37B10; Secondary: 37E05} \keywords{substitutive systems, ergodic sums, interval exchange transformations, Markov approximation}

\begin{abstract}
Deviation of ergodic sums is studied for substitution dynamical systems with a matrix that admits eigenvalues of modulus 1.  The functions $\gamma$ we consider are the corresponding eigenfunctions.
In Theorem \ref{eigenvalue1} we prove that the limit inferior of the ergodic sums
$(n, \gamma(x_0) + \ldots +\gamma(x_{n-1}))_{n\in \NN}$
is bounded for every point $x$ in the phase space.
In Theorem \ref{fluctuations-periodic}, we prove existence of limit distributions
along certain exponential subsequences of times for substitutions of constant length.
Under additional assumptions, we prove that ergodic integrals satisfy the Central Limit Theorem
(Theorem \ref{fluctuations-generic-basic}, Theorem \ref{fluctuations-generic}).
\end{abstract}

\date{14 June  2011}

\maketitle \markboth{ Xavier Bressaud, Alexander I. Bufetov, Pascal Hubert}{Substitutions with eigenvalues of modulus 1}
\tableofcontents

 \section{Introduction}

Let $\sigma$ be a primtive substitution over a finite alphabet $\A$, let $M_\sigma$ be the matrix substitution, and let $\xs$ be the corresponding subshift in the space of bi-infinite sequences (the precise definitions are recalled in Section  \ref{sec:background}).  The aim of this paper is to study the asymptotic behaviour of ergodic sums for the minimal dynamical system $(\xs, T)$  in the (non-hyperbolic) case when  the matrix $M_\sigma$ has an eigenvalue of modulus one.
For a function $f : X_\sigma \to \RR$, all $x \in X_\sigma$ and all $n \in \NN$, set
$$ S_n f (x)= \sum_{k=0}^{n-1} f (T^k x).$$
In this paper, we shall only consider functions $f$ depending on the first coordinate of the symbolic sequence. In what follows, we will identify such a function $f$ with the corresponding vector in $\CC^{\# \A}$.

\subsection{Limit inferior of ergodic sums}


\begin{theo} \label{eigenvalue1}
Let $\sigma$ be a primitive substitution with matrix $M_{\sigma}$.
 Assume that $M_{\sigma}$ has an eigenvalue $\theta_2$ of modulus 1 and let $\Gamma = (\gamma(a), a \in \A)^t$ be a right eigenvector associated to $\theta_2$. There exists a constant $C$ (depending only on $\sigma$ and $\gamma$) such that, for any $x \in \xs$, we have
\begin{equation} \label{bounded}
\liminf_{n\to \infty} \vert S_n \gamma (x) \vert < C
     \hbox{  and  }
\liminf_{n\to \infty} \vert S_n \gamma (T^{-n} x) \vert < C.
\end{equation}
\end{theo}

Recall that, since $T$ is minimal, by the Gottschalk-Hedlund Theorem, the function $\gamma$ is a coboundary if and only if the ergodic sums
$S_n \gamma (x)$
remain bounded, as $n\to\infty$, for some  and, consequently, for every $x\in \xs$.
Adamczewski \cite{Ad} proved that, if $\theta_{2}$ is not a root of unity, then  for almost all $x$, we have
\begin{equation}\label{eq:limsupinf}
\limsup_{n\to \infty} \vert S_n \gamma (x) \vert  = \infty.
\end{equation}
If $\theta_{2}$ is a root of unity, then Adamczewski \cite{Ad} gives a necessary and sufficient  combinatorial condition ensuring that $\gamma$ is a coboundary.
 For {\it almost all} points the conclusion of Theorem \ref{eigenvalue1}  is true for any ergodic transformation. It is important for us to have the result for all points.
 Observe that there exist infinitely many primitive substitutions with a hyperbolic matrix with the following property (see \cite{BHM}).
 If the function $\Gamma$ is in the unstable space of the matrix or equivalently if there exists $x$
 such that
  $$
 \limsup_{n\to \infty} \vert S_n \gamma (x) \vert = \infty
 $$

 then there exists $x$ such that

 $$
 \liminf_{n\to \infty} \vert S_n \gamma (x) \vert = \infty
     \hbox{  and  }
\liminf_{n\to \infty} \vert S_n \gamma (T^{-n} x) \vert = \infty
 $$
 Below, we give further applications of Theorem \ref{eigenvalue1} to affine interval exchange transformations (Corollaries \ref{AIET} and \ref{examplesIET}).

\subsection{Limit distributions for substitutions of constant length}

Informally, the mechanism of the Central Limit Theorem for substitutions with eigenvalue $1$ can be formulated as follows.
A long orbit of a substitution dynamical system admits a representation as a disjoint union of segments of different exponential scales.
Since the eigenvalue is equal to $1$, these segments give roughly the same contribution. By the Markov property,
these segments are asymptotically independent, whence the Central Limit theorem.

We are only able to make this argument precise under significant additional restrictions.

Recall that a substitution $\sigma$ has constant length $d$ if the image of  each letter has length $d$.  For instance, this condition  implies that the maximal eigenvalue of the matrix $M_\sigma$ is $d$. Furthermore, we assume that 1 is an eigenvalue and consider the eigenvector $\gamma$, $\gamma M_\sigma = \gamma$.  For a given $t \in [0,1]$, we consider ergodic sums $S_{\lfloor d^nt \rfloor }(\gamma)$ as a sequence of random variables where the initial point is distributed according to the unique shift invariant measure on $X_\sigma$. The behaviour of these sequences significantly depends on the expansion of $t$ in basis $d$;  namely on the sequence $\tau_t=(\tau_n)_{n \in \NN}  \in \{0, \ldots, d-1\}^\NN$ such that  $t= \sum_{n \leq 1}\tau_{n}d^n$. We first consider the case when this expansion is eventually periodic.
 \begin{theo} \label{fluctuations-periodic}
 Let $\sigma$ be a primitive substitution of constant length $d$. Assume that $M_\sigma$ has eigenvalue 1 with eigenvector $\gamma$.
 Let $t \in [0,1]$ and assume that the expansion $\tau_t =(\tau_{n})_{n \in \NN} $ of $t$ in basis $d$ is eventually periodic. Then,
 $$\frac{1}{\sqrt{n}} S_{\lfloor d^nt \rfloor }(\gamma) \xrightarrow{distr} \nu$$
 where $\nu$ is a probability measure with generalized density
 $$p_{0}\delta_{0} + \sum_{k = 1}^{\ell} \frac{p_{k}}{\sqrt{2 \pi } \sigma_{k}} e^{-x^2/(2 \sigma_{k}^2)}, \ \sum_{k = 0}^{\ell} p_{k} = 1.$$
 Moreover, the limit is trivial if and only if  $\gamma$ is a coboundary.
 \end{theo}
The nature of the limit law depends on the expansion of $t$ in basis $d$. In Section \ref{section:examples} we give examples showing that the limit distribution need not be Gaussian and can have atoms.

In Section  \ref{proof-thm-generic} we provide a precise algorithmically verifiable sufficient condition on the morphism $\sigma$ under which the limit distribution is indeed Gaussian for almost all $t$.
We thus arrive at the following theorem.




  \begin{theo}  \label{fluctuations-generic-basic}
 There exist infinitely many primitive substitutions $\sigma$ of constant length with eigenvalue $1$ and eigenvector $\gamma$ such that, for Lebesgue almost every $t \in [0,1]$, the following holds:
 \begin{enumerate}
 \item there exist positive constants $\alpha_{1}$, $\alpha_{2}$ such that  the variance  $V_{n} = Var(S_{\lfloor d^nt \rfloor }(\gamma))$ satisfies
 the inequality
 $$\alpha_1 n^{4/5} \leq V_n \leq \alpha_2 n. $$
 \item  The sequence $$\frac{1}{\sqrt{V_{n}}} S_{\lfloor d^nt \rfloor }(\gamma)$$
 converges in distribution to the normal law $N(0,1)$.
 \end{enumerate}
  \end{theo}

Our argument is an adaptation to our context of R. L. Dobrushin's central limit theorem \cite{Dob} for non homogeneous Markov chains.
We follow the martingale approach of S. Sethuraman and S.R.S. Varadhan \cite{SV}.

\begin{con}
$\mbox{}$
\begin{enumerate}
\item For any primitive substitution of constant length with eigenvalue 1,
there exists $t$ such that the limit distribution of $$\frac{1}{\sqrt{V_{n}}} S_{\lfloor d^nt \rfloor }(\gamma)$$ is not Gaussian.
\item For any primitive substitution of constant length with eigenvalue 1, for almost every $t$,  $$\frac{1}{\sqrt{V_{n}}} S_{\lfloor d^nt \rfloor }(\gamma)$$
 converges to the normal law $N(0,1)$.

\end{enumerate}

\end{con}

\subsubsection{Affine interval exchange transformations }
Our study was partly motivated by questions arising in the study of so-called Affine Interval Exchange Transformations. An interval exchange  transformation is defined by the length $\lambda = (\lambda_1, \ldots, \lambda_r)$ of the intervals $I_1, \ldots, I_r$ and a permutation $\pi$ of $\{1, \ldots, r\}$.  It is denoted by $T_{(\lambda, \pi)}$.
An affine interval exchange transformation is a piecewise affine bijective map of the interval. As for interval exchange transformation it is defined by the length of the intervals, a permutation, and the slopes of the map (constant on each interval). This additional information is encoded in a vector $(w_1, \ldots, w_r)$ with $w_i >0$ for $i = 1, \ldots, r$. Given these data, one can determine $\{a_i, i = 1, \ldots, r\}$ such that   for all $i =1, \ldots, r$ and all $x \in I_i$,
$$T(x) = w_i x + a_i.$$
To an affine interval exchange transformation $T$, we assign the corresponding $\log$-slope vector $$\vec{\Gamma} = (\log w_1, \ldots, \log w_r).$$
A natural idea to analyze this class of dynamical systems  has been to ask whether they are semi-conjugate to some interval exchange transformation (as a circle diffeomorphism is semi-conjugate to a rotation).  A subsequent question is to ask if such a semi-conjugacy is a conjugacy. Indeed, it is not always the case : affine interval exchange transformation may have wandering intervals.
It is well-known that the problem of conjugacy reduces to estimates on the ergodic sums of interval exchange transformations
(see for instance \cite{Ca}, \cite{Co}, \cite{LM}, \cite{BHM},  \cite{MMY2}). Together with this observation, our result yields new classes of examples corresponding to self-similar interval exchange transformation with eigenvalues of modulus one.

\vip

Let $T$ be an affine interval exchange transformation. One can assign to $T$ its {\it Rauzy expansion} in the same way as for usual interval exchange transformations \cite{MMY2}. Let us consider the case where the Rauzy expansion of $T$ is periodic with positive Rauzy renormalization matrix $R$. In this case $T$ admits a semi-conjugacy onto the underlying self-similar interval exchange transformation $T_{0}$. Let  $(\lambda_{1},  \dots \lambda_{n})$ be the lengths of $T_{0}$ and as before, let  $(w_{1}, \dots, w_{n})$ be the slopes of $T$. Camelier and Gutierrez  \cite{Ca} showed furthermore that
\begin{equation} \label{eq:orthog}
\sum_{i}^n \lambda_{i} \log w_{i} =  0.
\end{equation}
Conversely, given a self-similar interval exchange transformation $T_{0}$ and a vector of slopes $(w_{1}, \dots, w_{n})$ satisfying  \eqref{eq:orthog} then Camelier and Gutierrez showed that there exists an affine interval exchange transformation $T$ with slopes $(w_{1}, \dots, w_{n})$ semi-conjugate to $T_{0}$.

We next consider the question of conjugation.
Consider the decomposition of $\CC^n$ into characteristic spaces of the matrix $^t R$:
$$\CC^n = E^{u}\oplus E^c \oplus E^s.$$
Here $E^u$ corresponds to eigenvalues with absolute value greater than $1$, $E^s$  to eigenvalues with absolute value less than $1$,
$E^c$  to eigenvalues on the unit circle.
If $\vec{\Gamma} \in E^s$, then Camelier and Gutierrez \cite{Ca} showed that $T$ is conjugate to $T_{0}$. If $\vec{\Gamma} \in E^u$,
Camelier and Gutierrez \cite{Ca} and Cobo \cite{Co} gave examples in which $T$ has wandering intervals. In \cite{BHM} , it is showed that if $\vec{\Gamma} \in E^u$ and Galois conjugate to the Perron-Frobenius vector of $^t R$ then there is no conjugacy. Marmi-Moussa-Yoccoz \cite{MMY2} obtained the same conclusion for generic interval exchange transformation (\cite{MMY2}) using the positivity of the Lyapunov exponent of the Kontsevich-Zorich cocycle (proved by Forni \cite{Fo}, Avila-Viana \cite{AV}).

Theorem \ref{eigenvalue1} implies

\begin{cor} \label{AIET}
Let $T_{0}= T_{(\lambda, \pi)}$ be a self-similar  interval exchange transformation and $R$ the associated matrix obtained by Rauzy induction. Assume that $R$ has an eigenvalue $\theta$ of modulus 1.
\begin{itemize}
\item If $\theta = \pm 1$, let $\Gamma$ be an eigenvector associated to $\theta$ for $^t R$
\item If $\theta = e^{i\phi}$ belongs to $\CC \setminus \RR$, let $V_{\phi}$
(resp. $V_{-\phi}$) be an eigenvector associated to $\theta$ (resp. to $e^{-i\phi}$) for $^t R$. Let $\Gamma$ be a vector with real entries belonging to the vector space generated by $V_{\phi}$
 and $V_{-\phi}$
\end{itemize}
 If an  affine interval exchange transformation
with vector of the logarithm of the slopes is $\Gamma$ is semi-conjugate to $T_{0}$ then it is also  topologically conjugated to $T_{0}$.
\end{cor}

This corollary applies to many examples:

\begin{cor} \label{examplesIET}

There exist infinitely many self similar interval exchange transformation $T_{(\lambda, \pi)}$ on $r \geq 4$ intervals such that every affine interval exchange transformation that is semi-conjugate to $T_{(\lambda, \pi)}$ is conjugate to $T_{(\lambda, \pi)}$.
\end{cor}

We recall \cite{veechsurf} that the Veech group $\Gamma (X)$ of a translation surface $X$ is a discrete subgroup of $\textrm{SL}_{2}(\RR)$  obtained as the image by derivation of the group of affine diffeomorphisms of the translation surface $X$.
The surface $X$ is a Veech surface if its Veech group is a lattice in $\textrm{SL}_{2}(\RR)$. An element of $\Gamma$ is the derivative of a pseudo-Anosov diffeomorphism if and only if it is a hyperbolic element. A Veech surface is primitive if it is not a cover of any other translation surface. Assume that $X$ is a primitive Veech surface of genus 2. The group $\Gamma (X)$ is then
defined over a quadratic field $K$ (\cite{Mc}).  We have the following proposition:

\begin{prop} \label{prop:Veech-Salem}
On each primitive Veech surface in genus 2, there exists a pseudo-Anosov diffeomorphism whose dilatation has two conjugates of modulus one.
 \end{prop}
{\bf Remark.} In other words, the dilatation of our pseudo-Anosov diffeomorphism is  a Salem number.

 To a pseudo-Anosov diffeomorphism given by Proposition \ref{prop:Veech-Salem}, assign an interval exchange transformation using Veech's zippered rectangle construction. The resulting self-similar interval exchange transformation $T_{0}$ then has the desired property:
 any affine interval exchange transformation $T$ semi-conjugate to $T_{0}$ is in fact conjugate to $T_{0}$. Indeed the unstable space is spanned by the Perron-Frobenius vector and by the result of Camelier-Gutierrez (see equation \eqref{eq:orthog}), the vector of slopes cannot belong to the unstable space. If the vector of slopes belongs to the stable space, then again results of Camelier and Gutierrez imply the existence of a conjugacy.  If the vector of slopes belongs to the central space then Corollary \ref{AIET} applies and again yields the existence of a conjugacy.
Others examples can be obtained using the square tiled surfaces with degenerate  Lyapunov spectrum constructed by Forni \cite{Fo2} and Forni-Matheus-Zorich \cite{FMZ}.

\subsubsection{Symbolic flows and deviation of ergodic averages}

To prove Theorem \ref{fluctuations-periodic} we introduce a family of Markov chains which provide approximations of the ergodic sums. To give more specific statements ---but without entering now into details of the proof--- let us  consider, for each  $N \in \NN$ and $\tau_0\in \{0, \ldots, d^{N}\}$, the Markov chain $X^{N,\tau_0}$ defined on the state space $\A \times \A^2 \times \{0, \ldots, d^{N}\}$ having probability transitions:
$$p_{(a,V,m) \to (b,W,k)} = \left\{ \begin{array}{ll} d^{-N} & \hbox{ if } \sigma^N(a) = PbS \hbox{  with } |P|=m \\ & \hbox{ and } \sigma^N(V) = P'WS' \hbox{  with } |P'| = m+\tau_0\\ 0 & \hbox{ otherwise}. \end{array} \right.$$
Consider the map $G_{\tau_0}(a,V,m)= \gamma(S) + \gamma(P')$ ; we say that $G_{\tau_0}$ is a coboundary if it remains bounded
on all paths of the Markov chain $X^{N,\tau_0}$.   

 \begin{cor}
 \label{criterion}
Assume that the chain $X^{N,\tau_0}$ only possesses one recurrent component and set $t_0 = \sum_{n\geq1} \tau_0 d^{-nN}$. If $S_{n}\gamma$ is unbounded, then, as $n\to\infty$, the sequence of random variables $S_{\lfloor d^n t_0\rfloor}\gamma/\sqrt{n}$ converges in distribution to a normal law with positive variance.
 \end{cor}

 We prove by analyzing examples that the irreducibility assumption does not always hold and that the behaviour of the fluctuations
 at exponential scale can be very complicated. For instance, we exhibit examples of substitutions $\sigma$ with an invariant vector $\gamma$ that is not a coboundary such that $S_{d^nt }(\gamma)$ is bounded
for some  values of $t$.
 The combinatorial properties of the graphs of the corresponding Markov chains remain very unclear.
Nevertheless, using results on non stationary  Markov chains, we can prove that if  there is $t_0$ with a periodic expansion such that the associated Markov chain is aperiodic with positive variance, then fluctuations are normal for Lebesgue almost every $t \in [0,1]$.
More specifically, we prove the following Theorem.
\begin{theo}  \label{fluctuations-generic}
Let $\sigma$ be  a primitive substitution of constant length with eigenvalue $1$. Let $\gamma$ be an eigenvector associated with the eigenvalue $1$. Assume that there exists $N$ and $\tau_0$ such that the Markov chain $X^{N,\tau_0}$ is recurrent aperiodic and the function $G_{\tau_0}$ is not a coboundary.  Then, for Lebesgue almost every $t \in [0,1]$,
the following holds:
 \begin{enumerate}
 \item there exist positive constants $\alpha_{1}$, $\alpha_{2}$ such that  the variance  $V_{n} = Var(S_{\lfloor d^nt \rfloor }(\gamma))$ satisfies
 the inequality
 $$\alpha_1 n^{4/5} \leq V_n \leq \alpha_2 n. $$
 \item  The sequence $$\frac{1}{\sqrt{V_{n}}} S_{\lfloor d^nt \rfloor }(\gamma)$$
 converges in distribution to the normal law $N(0,1)$.
 \end{enumerate}

 \end{theo}
The sequence $V_{n}$ is the variance of a non stationary Markov chain.
The assumption (positivity of the variance for a specific $\tau_0$) means that for the homogeneous Markov chain associated to $t_{0}=\sum_{n\geq 1} \tau_0 d^{-nN}$, the approximation of the ergodic sums of $\gamma$ on this chain is not a coboundary. It is not very explicit but  it is possible to check it on examples.
 We prove that the hypothesis of Theorem \ref{fluctuations-generic} is satisfied in many cases.

\subsection{Historical remarks.}
For deviations of ergodic averages of generic interval exchange transformations  see \cite{Zo}, \cite{Fo}.
For holonomy flows of pseudo-Anosov diffeomorphisms for which the second eigenvalue $\theta_{2}$ of the corresponding action in homology is real and satisfies the inequality $\theta_{2} >1$, limit theorems are obtained in \cite{sasha2}. In this case, limit distributions have compact support.

Let $\sigma$ be a primitive substitution over a finite alphabet $\A$ and $\gamma$ be a function from $\A$ to $\CC$. Let $\xs \subseteq \A^\ZZ$ be the subshift defined from $\sigma$ and $x$ be a point in $\xs$. The stepped line associated to $x$ is the piecewise affine curve $\ell(x)$ in $\RR\times \CC$
 with vertices $(n, \gamma(x_0) + \ldots +\gamma(x_{n-1}))_{n\in \NN}$ and $(n, \gamma(x_{-n}) +\ldots + \gamma(x_{-1}))_{n\geq 1}$.
 Stepped lines were studied by many authors (see for instance \cite{ArIt} and \cite{ABB} and references there).
They are closely related to the famous Rauzy fractals
 and to the fractal curves studied by Dumont and Thomas in
\cite{DT1}, \cite{DT2}.
Using a different language, the result of Adamczewski \cite{Ad} about discrepancy of substitutive systems describes the behaviour of
\begin{equation}
\limsup_{n\to \infty} (\vert \gamma(x_0) +\ldots + \gamma(x_n)\vert)
\end{equation}
\noindent in terms of the eigenvalues of the matrix $M_\sigma$ of $\sigma$.

\subsection{Organization of the paper}
Section \ref{sec:background} sets our notation on words, sequences, substitutions and suspension flows.
Theorem \ref{eigenvalue1} is proved in Section \ref{proof-eigenvalue1} while its corollaries \ref{AIET} and \ref{examplesIET} are proved in Section \ref{proof-cor}. In Section \ref{section:automata}, we prove preliminary results about approximation of ergodic sums by Markov chains. More precisely, we construct an explicit family of finite automata which allow us to code the ergodic sums as sums of nonhomogeneous  Markov random variables. The main difficulty is that the underlying Markov chain may fail to be irreducible (examples are given in Section \ref{section:examples}).
In Section \ref{proof-fluctuations-periodic} we prove Theorem \ref{fluctuations-periodic}.  In Section \ref{proof-thm-generic} we prove Theorem \ref{fluctuations-generic}. In Section \ref{section:examples} we discuss various examples including those proving Theorem \ref{fluctuations-generic-basic} (see Remark \ref{remark-family}).

\subsection{Acknowledgements.}

We thank S. Geninska for mentioning Beardon's theorem \cite{Be}.
X.~B. was partially supported by  project  ANR JCJC: LAM while this work was in progress.
A.~I.~B. is an Alfred P. Sloan Research Fellow.
He is supported in part by~grant MK-4893.2010.1 of
the President of the Russian Federation,
by the Programme on Mathematical Control Theory of the Presidium of the Russian Academy of Sciences,
by the Programme 2.1.1/5328 of the Russian Ministry of Education and Research,
by the Edgar Odell Lovett Fund at Rice University,
by the NSF under grant DMS 0604386,
and by the RFBR-CNRS grant~\mbox{10-01-93115}.
P.~H. was partially supported by  project blanc ANR: ANR-06-BLAN-0038 while this work was in progress. 

\section{Background} \label{sec:background}
\subsection{Words and sequences}

Let $\A$ be a finite set. One calls it an {\it alphabet} and its elements {\it symbols}.
A {\it word} is a finite sequence of symbols in $\A$, $w=w_0\ldots w_{\ell-1}$. The length of $w$ is denoted $|w|=\ell$. One also defines the empty word $\varepsilon$.  The set of words in the alphabet $\A$ is denoted $\A^*$ and $\A^+=\A^*\setminus \{\varepsilon\}$.  We will need to consider words indexed
by integer numbers, that is, $w=w_{-m} \ldots w_{-1}.w_0\ldots w_{\ell}$ where $\ell,m \in \NN$ and the dot separates negative and non-negative coordinates.

The set of one-sided infinite sequences $x=(x_i)_{i \in \NN}$ in $\A$ is denoted by $\A^\NN$. Analogously, $\A^\ZZ$ is the set of two-sided infinite sequences
$x=(x_i)_{i\in \ZZ}$.

Given a sequence $x$ in $\A^+$, $\A^\NN$ or $\A^\ZZ$ one denotes $x[i,j]$ the subword of $x$ appearing between indexes $i$ and $j$. Similarly one defines $x(-\infty,i]$ and $x[i,\infty)$.
Let $w=w_{-m} \ldots w_{-1}.w_0\ldots w_{\ell}$ be a word on $\A$. One defines the cylinder
set $[w]$ as $\{ x \in \A^\ZZ : x[-m,\ell]=w \}$.

The shift map $T:\A^\ZZ \to \A^\ZZ$ or $T:\A^\NN \to \A^\NN$ is given by $T(x)=(x_{i+1})_{i \in \NN}$ for
$x=(x_i)_{i\in \NN}$. A subshift  is any shift invariant and closed (for the product topology) subset of $\A^\ZZ$ or $\A^\NN$. A subshift is minimal if all of its orbits by the shift are dense.

In what follows we will use the shift map in several contexts, often in restriction to a subshift.
To simplify notation we keep the symbol $T$ all the time.

\subsection{Substitutions} \label{subst} We refer to \cite{Qu} and \cite{fogg} and references therein for the general theory of substitutions.

A {\it substitution} is a map  $\sigma: \A \to \A^+$. It naturally extends to $\A^+$, $\A^\NN$ and $\A^\ZZ$; for
$x=(x_i)_{i\in \ZZ} \in \A^{\ZZ }$ the extension (which is a morphism of monoid) is given by $$
\sigma(x)=\ldots \sigma(x_{-2})\sigma(x_{-1}). \sigma(x_0)\sigma(x_1) \ldots
$$
where the central dot separates negative and non-negative coordinates
of $x$. A further natural  convention is that  the image of the empty word
$\varepsilon$ is $\varepsilon$.

Let $M$ be the matrix with indices in $\A$ such that $M_{ab}$ is the number of times letter $b$ appears in $\sigma(a)$ for any $a,b \in \A$. The substitution is {\em primitive} if there is $N >0$ such that for any $a \in \A$, $\sigma^N(a)$ contains any other letter of $\A$ (here $\sigma^N$ means $N$ consecutive iterations of $\sigma$). Under primitivity one can assume without loss of generality that $M>0$.

A substitution $\sub$ is of {\it  constant length} $\dd$ if, for all $a \in \A$, $|\sub(a)| = d$.

Let $\xs \subseteq \A^\ZZ$ be the subshift defined from $\sigma$. That is,
$x \in \xs$ if and only if any subword of $x$ is a subword of $\sigma^N(a)$ for some $N \in \NN$ and $a \in \A$.

Assume $\sigma$ is primitive. Given a point $x \in \xs$ there exists a unique sequence $(p_i,c_i,s_i)_{i\in \NN} \in (\A^* \times \A \times \A^*)^\NN$ such that for each $i \in \NN$: $\sigma(c_{i+1})=p_ic_is_i$ and
$$\ldots \sigma^3(p_3) \sigma^2(p_2)\sigma^1(p_1) p_0 . c_0 s_0 \sigma^1(s_1)\sigma^2(s_2)\sigma^3(s_3)\ldots$$
is the central part of $x$, where the dot separates negative and non-negative coordinates. This sequence is called the prefix-suffix decomposition of $x$ (see for instance \cite{CS}).

If only finitely many suffixes $s_i$ are nonempty, then there exists $a \in \A$ and non-negative integers $\ell$ and $q$ such that
$$x^+ = x[0,\infty)=c_0 s_0 \sigma^1(s_1)\ldots \sigma^\ell(s_\ell)  \lim_{n\to \infty} \sigma^{n q}(a)$$
Analogously,  if only finitely many $p_i$ are non empty, then
$$x^- = x(-\infty,-1]=( \lim_{n\to \infty} \sigma^{n p}(b) ) \sigma^m(p_m) \ldots \sigma^1(p_1) p_0$$
for some $b \in \A$ and non-negative integers $p$ and $m$.

\subsection{Vershik automorphisms and suspension flows} \label{subsection:Vershik}

In this subsection, we recall the construction of Vershik's automorphisms \cite{Ver}, \cite{VL} and their
continuous analogues \cite{Ito}, \cite{sasha}.
We refer to Section 2, 3, 5 of \cite{sasha} for details and further references. Given an oriented graph $\Gamma$ with $m$ vertices, let $\E (\Gamma)$ be the set of edges of $\Gamma$. For $e \in \E (\Gamma)$ we denote $I(e)$ its initial vertex and $F(e)$ its terminal vertex. To the graph $\Gamma$ we assign a non-negative $m\times m$ non-negative matrix $A(\Gamma)$ by the formula
\begin{equation}
\label{matrixa}
A = A(\Gamma)_{i,j} = \sharp\{e \in \E(\Gamma) \ : \ I(e) = i, F(e) = j\}
\end{equation}
We assume that $A$ is a primitive matrix (in other words, that the corresponding topological Markov chain is irreducible and aperiodic).

We define the Markov compactum:
$$Y = \{ y = y_{1} \dots y_{n}\dots : y_{n} \in \E(\Gamma), F(y_{n+1}) = I(y_{n})\}$$
The shift on $Y$ is denoted by $\frak{S}$.
Assume that there is an order on the set of edges starting from a given vertex. This partial order extends to a partial order on $Y$: we write $y<y'$ if there exists $l \in \NN$ such that $y_{l}< y'_{l}$ and $y_{n}= y'_{n}$ for $n >l$. The Vershik automorphism $T^Y$ is the map from $Y$ to itself defined by
$$T^Y y = \displaystyle\min_{y' >y}y'.$$
As $A$ is primitive, there is a unique probability measure invariant under $T^Y$ denoted by $\mu_{Y}$.

We now define a suspension flow over $(Y,T^Y)$ in the following way:
let $H$ be the Perron-Frobenius eigenvector of $A$. Then $h_{t}$ is the special flow over $(Y,T^Y)$ with roof function $\tau(y) = h_{I(y_{1})}$. The phase space of the flow is
$$Y(\tau) = \{(y,t) : y \in Y, Ê0 \leq t < \tau(y)\}.$$
The measure $\mu_{Y}$ induces a probability measure $\nu_{\Gamma}$ on $Y_{\tau}$. The vector $H$ is normalized in such a way that the
space $Y(\tau)$ have total measure $1$.

For each $e \in \E(\Gamma)$, the set $$ \{(y,t) : y \in Y, y_{1} = e, \ ÊÊ0 \leq t < h_{I(y_{1})}\}$$ is  called a {\it rectangle} in what follows.

The space $X =  \{ x = \dots x_{-n} \dots x_{n}\dots : x_{n} \in \E(\Gamma), F(x_{n+1}) = I(x_{n})\}$ is the natural extension of $(Y,S)$. The space $X$ endowed with the Parry measure (the measure of maximal entropy) is
canonically isomorphic to  $(Y(\tau), \nu_{\Gamma})$.

Following Livshits \cite{Liv},
we now connect Vershik's automorphisms with  substitution dynamical systems. Consider the alphabet $\A = \{1, \ldots, m\}$ as the set of vertices of $\Gamma$. For all $a \in \A$, we denote $\sigma(a)$ the sequence  $\{F(e) : I(e) = a\}$ ordered with the partial order on $\{e : I(e) = a\}$.  The dynamical  system $(\xs,T)$ is a topological factor of the Vershik automorphism $(Y,T^Y)$. The semi-conjugacy is given by the prefix-suffix decomposition. Almost every $u \in \xs$ can be written in the form
$$u = \cdots \sigma^n(P_n) \cdots \sigma(P_1) P_0 . a_0 S_0 \sigma(P_1) \cdots \sigma^n(P_n) \cdots, $$ where, for all $n \in \NN$, $\sigma(a_{n+1}) = P_n a_{n} S_n$. Thus,  it is clear that such point correspond to the unique path in $Y$ such that, for all $n \geq 0$, $a_n = F(y_{n+1})$ (and  $=I(y_n)$ for $n>0$)  and $y_{n+1}$ has exactly $|P_n|$ predecessors in the partial order around $I(y_{n+1})$. We recall that the semi-conjugacy is not always a topological conjugacy because there may be multiple writings but it is a measurable conjugacy.

When $\sigma$ has constant length, the roof function of the suspension flow is constant because $^t(1, \dots, 1)$ is an eigenvector for $A$ associated to the maximal eigenvalue. 

We observe that
the matrix $M_\sigma$ of the substitution is connected to the matrix $A$ by the relation
$M_\sigma = A^T$.

\section{Proof of Theorem \ref{eigenvalue1}}
\label{proof-eigenvalue1}
Without loss of generality, we assume that $M$ is a positive matrix.
Let $x \in \xs$, $\theta$ be an eigenvalue of $M$ of modulus one and
 $\Gamma = (\gamma(a), a \in \A)$ an eigenvector associated to $\theta$.
We give a proof for
$$\displaystyle\liminf_{n\to \infty} \vert (\gamma(x_0) +\ldots + \gamma(x_n))\vert < C.$$

The other estimate is quite similar to obtain.

The prefix-suffix sequence $(p_i,c_i,s_i)_{i\in \NN}$ of $x$ is a path in a {\it finite} automaton (the prefix-suffix automaton, see \cite{CS}). Therefore, the prefixes and suffixes belong to a finite set $\P\times \S$ depending only on $\sigma$.

\medskip

{\bf First Case} We first assume that infinitely many suffixes are non empty in the prefix-suffix decomposition of $x$.  Let $k$ be a positive integer, by the prefix-suffix decomposition, we have
$$\begin{matrix}
x = \ldots \sigma^{k+1}(p_{k+1})\sigma^n(p_n) \sigma^{k-1}(p_{k-1}) \dots \sigma^1(p_1) p_0 . c_0 s_0 \sigma^1(s_1)
 \\
\dots \sigma^{k-1}(s_{k-1})\sigma^k(s_k)\sigma^{k+1}(s_{k+1})\ldots
\end{matrix}$$

We denote by $\hat{P}_k$ the word $\sigma^{k-1}(p_{k-1}) \dots \sigma^1(p_1) p_0$ and by $\hat{S}_k$ the word $c_0 s_0 \sigma^1(s_1)\dots \sigma^{k-1}(s_{k-1})$.
Thus we have
$$
x = \ldots \sigma^{k+1}(p_{k+1})\sigma^k(p_k) \hat{P}_k . \hat{S}_k\sigma^k(s_k)\sigma^{k+1}(s_{k+1})\ldots$$
\noindent and
$$ \sigma^k({c_k}) = \hat{P}_k . \hat{S}_k.$$

Assume that $s_{k+1}$ is non empty. Then $\sigma(s_{k+1})$ contains $c_k$ because every letter appears in the image of every letter. We consider the first occurrence of $c_k$ in $\sigma(s_{k+1})$ and get
$$\sigma(s_{k+1}) = \Pi_k c_k \Sigma_k,$$
\noindent where $\Pi_k$ and $\Sigma_k$ are words of bounded length (the bound only depends on $\sigma$).
Consequently $\hat{S}_k\sigma^k(s_k)\sigma^k(\Pi_k)\sigma^k(c_k)\sigma^k(\Sigma_k)$
is a prefix of $x^+ = x_0 x_1 \dots x_n \dots$.
This word is equal to $\hat{S}_k\sigma^k(s_k)\sigma^k(\Pi_k)\hat{P}_k\hat{S}_k\sigma^k(\Sigma_k).$
Thus, the word $W_k = \hat{S}_k\sigma^k(s_k)\sigma^k(\Pi_k)\hat{P}_k$ is a prefix of $x^+$.

We extend the function $\gamma$ to the monoid $\A^*$ by the natural manner.
$$\gamma: \begin{array}{lcl}
\A^* & \to & \CC  \\
w= w_1 \dots w_t & \mapsto & \displaystyle\sum_{i=1}^t\gamma(w_i).
\end{array}
$$
It is clear that  the extension of $\gamma$ is a morphism of monoid: if $v$ and $w$ are two words $\gamma(vw)= \gamma(v) + \gamma(w)$.

The key lemma is the following
\begin{lem} \label{W_k}
The sequence $(\gamma(W_k))_{k\in \NN}$ is bounded by a constant that only depends on $\sigma$ and $\gamma$.
\end{lem}

\begin{proof}
By the hypothesis on $\theta$ and $\gamma$, for every $w \in \A^*$ and $m \in \NN$,
$$\gamma(\sigma^m(w)) = \theta^m \gamma(w).$$
This implies that $$\vert\gamma(\sigma^m(w))\vert = \vert \gamma(w)\vert.$$
We use this remark to calculate
$\gamma(W_k) = \gamma(\hat{S}_k) + \theta^k\gamma(s_k) + \theta^k\gamma(\Pi_k) + \gamma(\hat{P}_k) = \gamma(\hat{S}_k\hat{P}_k) + \theta^k\gamma(s_k) + \theta^k \gamma(\Pi_k).$

Thus, $\gamma(W_k) = \gamma(\hat{S}_k\Pi_k) + \theta^k(\gamma(s_k) + \gamma(\Pi_k)) =
\gamma(\sigma^k(c_k)) + \theta^k(\gamma(s_k) + \gamma(\Pi_k)) =
\theta^k(\gamma(c_k) + \gamma(s_k) + \gamma(\Pi_k)).$

This yields that
$$\vert\gamma(W_k)\vert \leq \vert \gamma(c_k) \vert + \vert \gamma(s_k) \vert + \vert \gamma(\Pi_k) \vert.$$

The length of the words $c_k$, $s_k$, $\Pi_k$ are bounded independently on $x$. Thus, each factor on the right hand side of the inequality is bounded and the sequence
$(\vert\gamma(W_k)\vert)_{k \in \NN}$ is bounded independently on $x$.

\end{proof}

Now we end the proof of Theorem \ref{eigenvalue1} in the first case. We already proved that  $(W_k)_{k \in \NN}$ is a sequence of prefixes of $x^+$. By hypothesis, the length of $W_k$ tends to infinity (because infinitely many  suffixes $s_n$ are non empty). By Lemma \ref{W_k} the sequence $(\vert\gamma(W_k)\vert)_{k \in \NN}$ is bounded independently on $x$. This proves
Theorem \ref{eigenvalue1} in the first case.

\bigskip

{\bf Second Case} If the suffixes of $x$ are eventually empty, it means that $x^+$ belongs to the negative orbit of a periodic point for $\sigma$.
The proof is similar to the previous one.
There exists $a \in \A$ and non-negative integers $\ell$ and $q$ such that
$$x^+ = c_0 s_0 \sigma^1(s_1)\ldots \sigma^\ell(s_\ell)  \lim_{n\to \infty} \sigma^{n q}(a).$$
Let $S_{\ell} = c_0 s_0 \sigma^1(s_1)\ldots \sigma^\ell(s_\ell) $; by construction of the prefix-suffix automaton,
there exists a letter $c$ and a finite word $P_{\ell}$ such that
$\sigma^{\ell+1}(c) = P_{\ell}S_{\ell}$, thus
$$x = \dots P_{\ell}.S_{\ell}  \lim_{n\to \infty} \sigma^{n q}(a).$$
Since the letter $c$ appears in the word $\sigma(a)$, there exists words of bounded length $\Pi$ and $\Sigma$ such that
$\sigma(a) = \Pi c \Sigma$.
Thus for each $n > \ell +2$,  $S_{\ell} \sigma^{n q}(a) = S_{\ell} \sigma^{n q - \ell - 2} \sigma^{\ell +1}(\Pi c \Sigma)$ is a prefix of $x^+$.
Moreover
$$ S_{\ell} \sigma^{n q - \ell - 2}( \sigma^{\ell +1}(\Pi )) \sigma^{n q - \ell - 2} (\sigma^{\ell + 1}( c)) =
S_{\ell} \sigma^{n q - \ell - 2} (\sigma^{\ell + 1}(\Pi )) \sigma^{n q - \ell - 2}(P_{\ell}S_{\ell}).$$

Set $W_{n} = S_{\ell} \sigma^{n q - \ell - 2} (\sigma^{\ell + 1}(\Pi )) \sigma^{n q - \ell - 2}(P_{\ell})$; it is a prefix of $x^+$.

Now  $\gamma(W_{n}) = \gamma(S_{\ell}) + \gamma(\Pi) + \gamma(P_{\ell}) = \gamma(\sigma^{\ell+1}(c)) +  \gamma(\Pi) =  \gamma(c) + \gamma(\Pi)$.

The absolute value of this quantity is bounded. The bound only depends on $\sigma$ and $\gamma$ since $c$ is a letter and $\Pi$ belongs to a finite set (prefixes of $\sigma(a)$).

\section{Proof of Corollaries \ref{AIET} and \ref{examplesIET}}
\label{proof-cor}
\subsection{ Proof of Corollary \ref{AIET}}

Let $T = T_{(\lambda, \pi)}$ be a self-similar interval exchange transformation on $r$ intervals and $\Gamma = (\gamma_1, \dots, \gamma_r)$ a vector orthogonal to $\lambda$. The set of affine interval exchange transformations semi-conjugated to $T_{(\lambda, \pi)}$ with slopes $w = (e^\gamma_1, \dots, e^\gamma_r)$ is never empty (see \cite{Ca}, \cite{Co}); this set is denoted by $S_{aff}(T,w)$.

Cobo gave an if and only if condition ensuring the existence of $f$ in $S_{aff}(T,w)$ which is not conjugated to $T$ (see \cite{Co} page 392--393). This condition is the following:

 there exists a point $x \in \xs$ such that
\begin{equation} \label{finitesum}
\sum_{n\geq 1} e^{{\gamma(x_{-n}\ldots x_{-1})}} + \sum_{n\geq 1} e^{{-\gamma(x_0\ldots x_{n-1})}} <\infty.
\end{equation}
Otherwise every element of $S_{aff}(T,w)$ is conjugated to $T$ by a homeomorphism.

Theorem \ref{eigenvalue1} says that, if $\Gamma$ is an eigenvector associated to an eigenvalue of modulus one, then for every $x \in \xs$ the sequence $(\gamma(x_0\ldots x_{n-1}))_{n \in \NN}$ possesses a bounded subsequence.
  Therefore, if $ \theta = \pm 1$, $$\sum_{n\geq 1} e^{{\gamma(x_{-n}\ldots x_{-1})}} + \sum_{n\geq 1} e^{{-\gamma(x_0\ldots x_{n-1})}}$$   diverges for every x $\in \xs.$

 If $\theta \in \CC \setminus \RR$, by assumption, there exist complex numbers $a$ and $b$ such that $\Gamma = a V_{\phi} + b V_{-\phi}$. We apply Theorem \ref{eigenvalue1} to  the sequences  $(V_{\phi}(x_0\ldots x_{n-1}))_{n \in \NN}$ and $(V_{-\phi}(x_0\ldots x_{n-1}))_{n \in \NN}$. It is clear in the proof of Theorem \ref{eigenvalue1} that one can find a subsequence $(n_k)$ such that  $(V_{\phi}(x_0\ldots x_{n_k-1}))_{k \in \NN}$ and $(V_{-\phi}(x_0\ldots x_{n_k-1}))_{k \in \NN}$ are simultaneously bounded.
 Thus equation \eqref{finitesum} cannot be fulfilled which completes the proof of Corollary
 \ref{AIET}.

\subsection{Proof of Corollary \ref{examplesIET}}

\begin{figure}
 \includegraphics[width=11cm]{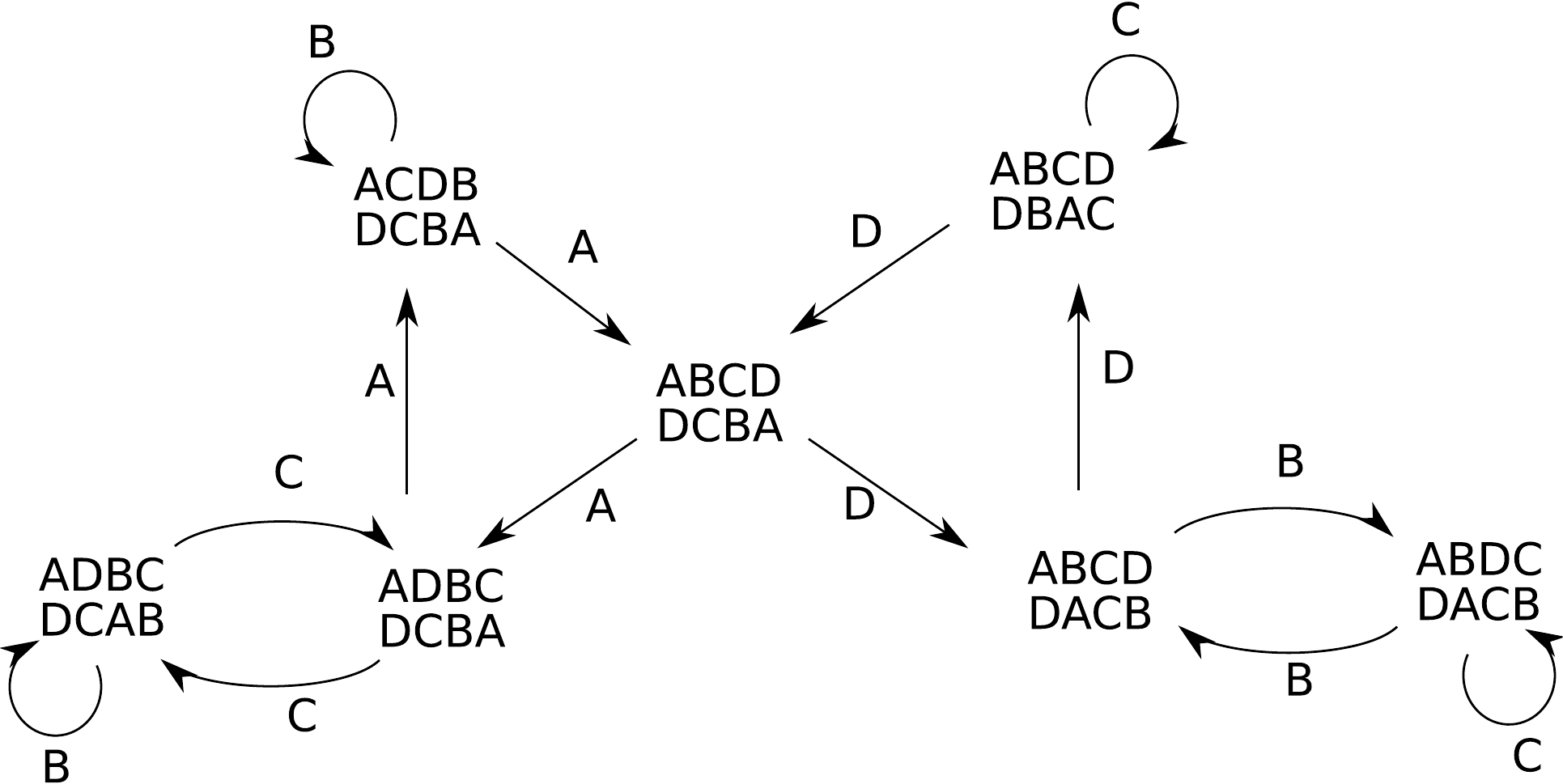}
\caption{\label{fig:rauzyclass}
Rauzy diagram.
}
\end{figure}

From Theorem \ref{eigenvalue1}, it is enough to
 give an infinite family of self-similar interval exchanges transformations on 4 intervals that have eigenvalues of modulus 1. The examples are obtained in the most simple Rauzy class (see figure \ref{fig:rauzyclass}).

We recall that a Rauzy class is the closure of a given pair of permutations under Rauzy induction and that self similar interval exchange transformations correspond to  loops in the Rauzy diagram. Marmi-Moussa-Yoccoz gave an if and only if condition ensuring that a loop is realized by a self similar interval exchange transformation. In our situation, the Marmi-Moussa-Yoccoz criterion says that a loop is admissible if all the labels $A$, $B$, $C$, $D$ appear along the loop (see \cite{MMY1}).
This condition will be fulfilled by the family of examples studied here.

Every arrow of the diagram induces a substitution. Given an interval exchange transformation $T$ on an interval $I$, let $T'$ be its image under Rauzy induction acting on $I'$.  The substitution is defined by following the images of elements of the partition $I'$ in the partition $I$ under $T$ until they come back to $I'$.
The product of these substitutions gives the substitution corresponding to the loop or to the self similar interval exchange transformation associated to the loop. Thus a loop is labelled by a finite word on the alphabet $\{A,B,C,D\}$.

The symmetric permutation is
$$\begin{pmatrix}
A&B&C&D \\
D&C&B&A
\end{pmatrix}$$

\begin{prop} \label{examples-EI-VP1}
For every $n \geq 1$, the matrix $M_{n}$ of the path starting from the symmetric permutation labelled by the loop $DBBDC^nDAAA$ has the following properties:

\begin{enumerate}
\item The Perron-Frobenius eigenvalue $\theta^{(n)}_{1}$ is an algebraic integer of degree 4.
\item It has two conjugate of modulus 1, $\theta^{(n)}_{2} = \overline{\theta^{(n)}_{3}}$
(it is a Salem number).
\item If $\theta^{(n)}_{2} = e^{2i\pi \alpha_{n}}$, then $\alpha_{n}$ is an irrational number.
\end{enumerate}
\end{prop}

\begin{proof}

The substitution on the alphabet $\{1,2,3,4\}$ is
$$\sigma_n = \left\{\begin{array}{ccl}
1 & \mapsto & 14 \\
2 &\mapsto & 14224 \\
3 &\mapsto & 14(23)^{n+1}24 \\
4 &\mapsto & 14(23)^{n}24.
\end{array}\right.
$$

The corresponding matrix is
$$M_n = \begin{pmatrix}
1 & 1 & 1 & 1 \\
0 & 2 & n+2 & n+1 \\
0 & 0 & n+1 & n \\
1 & 2 & 2 & 2
\end{pmatrix}.$$

The characteristic polynomial of $M_{n}$ is
$$P_{n} = X^4-(6+n)X^3 + (10+n)X^2-(6+n)X+1.$$
 By the general theory of interval exchange transformations (see \cite{Yo} \cite{Zo2}), $P_{n}$ is a reciprocal polynomial. Thus, the roots $\theta_{1}, \theta_{2}, \theta_{3}, \theta_{4}$ of $P_{n}$ satisfy
$\theta_{4} =  1/\theta_{1}$, $\theta_{3} =  1/\theta_{2}$ and $\theta_{1} > \vert \theta_{2}\vert \geq \vert \theta_{2}\vert > \theta_{4}$ (as $n$ is given, the dependence on $n$ is omitted in the sequel).

We first prove that $P_{n}$ has two complex roots of modulus 1. Let $t = \theta_{1} + 1/\theta_{1}$ and
$s = \theta_{2}+1/\theta_{2}$. We remark that $\theta_{2}$ and $\theta_{3}$ are complex conjugate if and only if $-2 \leq s \leq 2$. In fact, if $\theta_{2}$ is real and larger than 1, $s > 2$; if $\theta_{2}$ has modulus 1,
$\theta_{2} = e^{2i\pi \alpha}$ and $s = 2\cos (2\pi \alpha)$.

A simple calculation yields the equations $t+s = 6+n$ and $ts = 8+n$.
 Thus, $s = \frac{6+n -\sqrt{n^2 +8n+4}}{2}$.
We easily check that this quantity belongs to $(-2,2)$ independently on $n$. Thus $\theta_{2}$ has modulus 1.

Now, we prove that $P_{n}$ is irreducible. Otherwise, $\theta_{2}$ is $\pm 1$  or an algebraic integer of degree 2. A direct computation proves that $\pm 1$ are  not roots of $P_{n}$. Moreover, the only irreducible polynomials of degree 2, with roots of modulus 1 are the cyclotomic polynomials
$X^2+X+1$, $X^2+1$, $X^2-X+1$ (the polynomial has the form $X^2-sX+1$ with $-2<s<2$). To prove that $P_{n}$ is irreducible, we check that it is not a multiple of $X^2+X+1$, $X^2+1$, $X^2-X+1$.

We end the proof by showing that $\alpha$ is irrational.
If $\alpha$ is rational it means that $\theta_{2}$ is a root of unity. This is impossible since $\theta_{1}>1$ is a conjugate of $\theta_{2}$.

\end{proof}

\subsection{Proof of  Proposition \ref{prop:Veech-Salem}}
Denote by $\tau$ the non-trivial Galois automorphism of $K$ and by $\Gamma (X)$ the Veech group of $X$.
 Given
$A = \begin{pmatrix}
a & b \\
c & d
\end{pmatrix}
  \in \Gamma$, its Galois conjugate is $\tau(A) = \begin{pmatrix}
\tau (a) & \tau(b) \\
\tau(c) & \tau(d)
\end{pmatrix}$.
Denote by $\Gamma'$ the image of $\Gamma (X)$ under Galois conjugacy. Let  $\phi$ be an affine diffeomorphism of $X$ and $A$ the corresponding element of the Veech group,  its action on $H_{1}(X, \RR)$ is by
$\begin{pmatrix}
A & 0 \\
0 & \tau(A)
\end{pmatrix}$ in a suitable basis.

 The group $\Gamma'$ is a non discrete and non elementary subgroup of $\textrm{SL}_{2}(\RR)$.  By a result of Beardon \cite{Be}, it contains an elliptic element of infinite order. Let  $A^{\prime}=\tau(A^{\prime\prime})$ be this element. 
 If  $A^{\prime\prime}$ is hyperbolic, then  Proposition
 \ref{prop:Veech-Salem} follows. We now show that $A^{\prime\prime}$ cannot be elliptic nor parabolic.
 As $\Gamma (X)$ is discrete, every elliptic element of $\Gamma (X)$  is of finite order. But $A^{\prime}=\tau(A^{\prime\prime})$ is an element of infinite order. Consequently, $A^{\prime\prime}$  is not an elliptic element. Now the image of a parabolic element under Galois conjugacy  is again parabolic (since an element of $\textrm{SL}_{2}(\RR)$ is parabolic if and only if the absolute value of its  trace is 2), so $A^{\prime\prime}$ cannot be parabolic. Therefore $A^{\prime\prime}$ is indeed hyperbolic and is the derivative of a pseudo-Anosov diffeomorphism with the desired properties.

\section{Asymptotic behaviour of ergodic sums} \label{section:automata}

\subsection{Ergodic means}

We are interested in the distribution of the ergodic integral
$$\int_{0}^{\lambda_1^n t} f \circ h_s(x)  ds$$
considered as a random variable on the probability space $(X, \nu_{\Gamma})$ and normalized to have variance $1$.

The first observation is that for all $n \in \ZZ$, the law of $x$ is the same as the law of $\frak{S}^n(x)$ by $\frak{S}$-invariance of the measure $\nu_\Gamma$ (we mean: the law of $Id$ is the same as the law of $\frak{S}^n$). Hence the distribution of the ergodic integral 
is the same as that of the sequence 
$$\int_{0}^{\lambda_1^n t} f \circ h_s(\frak{S}^{-n} (x))  ds. $$
Let us state this result more formally:
\begin{lem} \label{invariance-under-frakS}
For any measurable function  $\varphi : Y \to \RR$,  all $z \in \RR$ and all $m \in \ZZ$,
$$ \nu( \{x \in Y : \varphi(x) \leq  z \} ) =  \nu( \{x \in Y : \varphi(\frak{S}^m x) \leq z \} ). $$
\end{lem}

From here till the end of the paper $\sigma$ is a substitution of constant length $d$, thus $\lambda_{1} = d$.
We choose a function $f$ constant on rectangles (defined in Subsection \ref{subsection:Vershik}). We are interested in the case when this function ``corresponds'' to the coordinates of an  eigenvector associated with eigenvalue $1$, i.e. an invariant vector (ie $f(\sigma(a)) = f(a)$ for every letter $a$).

\subsection{Decomposition of unstable leaves.}
Pieces of unstable leaves on which we compute the ergodic averages cross
transversally the ``rectangles'' of the partition.
 To such a piece of unstable leave we  associate a symbolic sequence given by the names of the successive rectangles it does cross. The self similar structure implies that these symbolic sequences are words in the language of the substitution $\sigma$. More precisely, such a piece of leaf starts inside a rectangle (the initial rectangle) then crosses transversally a given number of rectangles, to end up inside a last rectangle (the final rectangle).  Here, we try to describe the symbolic sequence using the self-similar structure. We have to take care of the initial and final rectangles which are not completely crossed. For (important) technical reasons we distinguish not only the final but the last two rectangles.

\vip

Let $x \in Y$ and $t \in [1; \infty)$. We consider a piece $\gamma_t^x$ of unstable leave determined by its starting point $x$ and its length $t$. The sequence $(x_n)_{n \in \ZZ}$ determines a bi-infinite sequence of prefixes/center/suffixes $(P^x_n,c^x_n,S^x_n)$ satisfying for all $n \in \ZZ$, $F(x_n) = c^x_n$ and
$\sigma(c^x_{n+1}) = P_n^x c_n^x S_n^x$. We observe that around $x$ the symbolic sequence (sequence of codes of rectangles seen along the unstable leave through $x$) is (note that this decomposition stops if all suffixes are empty)
$$\omega^x = \cdots \sigma^n(P^x_n) \cdots \sigma(P^x_1) P^x_0. c^x_0 S^x_0 \sigma(S^x_1) \cdots \sigma^n(S^x_n) \cdots.$$
We compute the distance from $x$ to the boundary of the rectangle in which it is by $t^-=\sum_{n <  0} \dd^n  |S^x_{n}|$.

We decompose (in basis $d$) : $t = \sum_{n \leq M} \dd^n \tau_n$ and let $k_t = \lfloor t \rfloor - 1 = \sum_{n=0}^{M} \dd^n \tau_n -1$.
We  observe that $\gamma_t(x)$ crosses at least $k_t$ rectangles. It does cross one more only if
$\sum_{n<0}  \dd^n  |P^x_{n}| +  \sum_{n <0} \dd^n \tau_n > 1$ (which happens only if
there is $H$ such that for all $h<H$, $|P^x_{-h}| + \tau_{-h} = d-1$, and $|P^x_{-H}| + \tau_{-H} \geq d$).

We denote $K=K_t$ the greater integer such that $d^K  < k_t$. We observe that $\ell_x = |S^x_0 \sub(S^x_1) \cdots \sub^K(S^x_K)| < k_t$. We also observe that there is a letter $W$ such that $|S^x_0 \sub(S^x_1) \cdots \sub^K(S^x_K) \sub ^{K+1}(W)| > k_t$ and that $S^x_0 \sub(S^x_1) \cdots \sub^K(S^x_K) \sub^{K+1}(W)$ is a factor of $\omega^x$. We decompose the prefix of length $k_t-\ell_x$ of  $\sub^{K+1}(W)$ into prefixes in such way that
$\sum_{i=0}^K d^i |S^x_i|  + \sum_{i=0}^K d^i |P^{t}_i|=k$ and $\prod_{i=0}^K \sub^i(S^x_i)  \prod_{i=K}^0 \sub^i(P^{t}_i) = \omega^x_1 \cdots \omega^x_{k_t}$.

Then, we set $t^+ = t - k - t^-$. We have to check if whether  $t ^+<1$ or $1\leq t^+<2$.
We use the identity $f(\sub(a)) = f(a)$ to get, if  $t ^+<1$,  the following equation:
\begin{equation} \label{equation:ergo-sum1}
\int_{0}^{t} f \circ h_s(x) ds = f(a_0) t^- + \sum_{n=1}^K f(S^x_n) + \sum_{n=K}^1 f(P^{t}_n) + t^+ f(a_k), 
\end{equation}
and, if $1\leq t^+<2$,
\begin{equation}  \label{equation:ergo-sum2}
\int_{0}^{t} f \circ h_s(x) ds = f(a_0) t^- + \sum_{n=1}^K f(S^x_n) + \sum_{n=K}^1 f(P^{x,t}_n) +  f(a_k) + (t^+-1)f(a_{k+1}).
\end{equation}
It is important to notice that except for the boundary terms, the decomposition depends on $x$ through $(x_n)_{n \geq 0}$ and of $t$ through $ \lfloor t \rfloor $. The other data tell about what happens more precisely in the boundary terms. We get
\begin{equation}
\left| \int_{0}^{t} f \circ h_s(x) ds - \left( \sum_{n=1}^K f(S^x_n) + \sum_{n=K}^1 f(P^{t}_n) \right) \right| \leq 3 ||f||_\infty. 
\end{equation}
If we step to $\dd^n t$ and shift the starting point to $\sigma^{-n}(x)$, we will obtain the same decomposition and something more precise on boundary terms. Hence, to get recursively to the next step we only have to keep track of $c_0^x$ and the word $a_{k_t}a_{k_t  +1}$. To encode this information, we define a family of automata which is described in the next paragraph.

\subsection{A family of automata}
The following family of automata associated with the substitution $\sigma$ will allow us
approximate the ergodic averages in a Markovian way.

\vip
 For all $W \in \A$ and all integer $1 \leq m \leq \dd$, we decompose  $\sigma(W) = P^{W,m} c^{W,m} S^{W,m}$ with $P^{W,m} \in \A^{m-1}$, $c^{W,m} \in \A$ and $S^{W,m} \in \A^{ \dd -m}$.

For all $\tau \in \{0, \ldots, \dd-1\}$ we define an automaton $\A_\tau$. States are couples $(a,W) \in \A \times \A^2$ and edges are labelled by a real number $v$ and an integer  $m$. For all $v \in \RR$, $m \in \{1, \ldots, \dd\}$  and integers and $a,b,  \in \A$, $V,W \in \A^2$ we put a labelled arrow  $$(a,V) \stackrel{v,m}{\to} (b,W)$$ if $b = c^{a,m}$, $W =  c^{V,m+\tau} c^{V,m+\tau+1}$ and  $v = \projinv(e_{S^{a,m}}+ e_{P^{V,m+\tau}}).$ Here the word $\sigma(V)$ has length $2d$, thus it has a prefix of length $m+\tau$. Observe that in particular,
$\sigma(a) = P^{a,m} b S^{a,m}$ so that  the outgoing degree $\dd$  of any vertex is the same as the outgoing degree of the corresponding vertex in the prefix/suffix automaton.  Note that for $\tau=0$ we can forget the second letter of $V$ and $W$;  indeed since $m\leq d$, $c^{V,m}$ is determined by the first letter of $V.$

Examples of such automata are described in Section  \ref{section:examples}. The readers should refer to this section to understand concrete examples.

\subsection{Ergodic sums and automata.}
We fix $t \in [1; \dd)$ and write $t = \sum_{n \leq 0} \tau_n \dd^n$ (note that $\tau_0 >0$). We let $W$ be a word of length $\dd$ i.e. $W \in \A^\dd$ and choose $x \in Y$ such that $\omega^x \in [W]$, i.e. $\omega^x_0 \cdots \omega^x_{\dd-1} = W$ or again $F(x_0)F((T_Yx)_0) \cdots F((T_Y^{\dd-1}x)_0) = W$.   The sequence $(x_n)_{n \in \ZZ}$ determines a bi-infinite sequence of prefixes/centres/suffixes $(P^x_n,c^x_n,S^x_n)$ satisfying for all $n \in \ZZ$, $F(x_n) = c^x_n$ and $\sub(c^x_{n+1}) = P_n^x c_n^x S_n^x$.

We observe that around $x$ the symbolic sequence (sequence of codes of rectangles seen along the unstable leave through $x$) is (except if all prefixes are empty)
$$\omega^x = \cdots \sub^n(P^x_n) \cdots \sub(P^x_1) P^x_0. c^x_0 S^x_0 \sub(S^x_1) \cdots \sub ^n(S^x_n) \cdots.$$
We note that $c^x_0=W_0$. In the simpler case, $W$ is the prefix of length $\dd$ of $c_0^x S^x_0 \sub(S^x_1)$ ; or of $c_0^x S^x_0 \sub^l(S^x_l)$ where $S^x_l$ is the first non empty suffix after $S_0^x$. If all suffixes are empty, then we must look at $T^k_Y(x)$. Anyhow, we set $a^0 := c^x_0 = W_0$ and $V^0  := W_{\tau_0} W_{\tau_0+1}$. We also set $U^0 = W_1 \cdots W_{\tau_0-1}$ (may be empty, if $\tau_0=1$).

\vip

We construct recursively a sequence $(a^n, V^n, U^n)_{n\geq 0}$. Assume $(a^n, V^n, U^n)$ determined for some $n$. Let then $m$ be the length $m = |P_{-(n+1)}^x|+1 (= \dd - |S^x_{-(n+1)}|)$.   We set
$$\left\{ \begin{array}{lcl}
a^{n+1} &:= & c^x_{-(n+1)} = c^{a^n,m} \\
V^{n+1} &:= &V^n_{\tau_{-n}+m} V^n_{\tau_{-n}+m+1} = c^{V^n, \tau_{-n}+m} c^{V^n, \tau_{-n}+m+1} \\\
U^{n+1} &:= & S_{-(n+1)}^x \sub(U^n) V^n_1 \cdots V^n_{\tau_{-n}+m-1}.
\end{array}
\right.$$

\begin{rem} This sequence $(a^n, V^n)$ is obtained by the family of automata described in the previous section.  The sequence of automata depends on the digits of $t$.
\end{rem}

Since $f$ depends only on the first coordinate, for brevity we write  $f(W)$ for the sum $f(W_1)+\cdots + f(W_k)$ if $W=W_1\cdots W_k$. 
We set $t_-^n = \dd^n \sum_{k<-n} \dd^k |S^x_{-k}|$ and $t_+^n =  \dd^n \sum_{k<-n} \dd^k (\tau_k +|P^x_{-k}|)$ and observe that $0 < t^n_+ < 2$. \begin{lem} \label{lem:link-sums-automaton}
If $ t_+^n < 1$, then 
$$ \int_{0}^{\dd^nt} f \circ h_s(\frak{S}^{-n} x) ds = f(a_n) t^n_- + f(U^n) + t^n_+ f(V^n_{1}).$$
If $ t_+^n \geq 1$, then
$$ \int_{0}^{\dd^nt} f \circ h_s(\frak{S}^{-n} x) ds = f(a_n) t^n_- + f(U^n) + f(V^n_{1}) + (t^n_+-1)  f(V^n_{2}).$$
In particular, we have 
$$ \left|  \int_{0}^{\dd^nt} f \circ h_s(\frak{S}^{-n} x) ds - f(U^n) \right|  \leq 3 ||f||_\infty.$$
\end{lem}

\begin{proof}
The symbolic decomposition around $\frak{S}^{m}(x)$ 
is
$$\omega^{\frak{S}^{m}x} = \cdots \sub^n(P^x_{m+n}) \cdots \sub(P^x_{m+1}) P^x_m. c^x_{m} S^x_m \sub(S^x_{m+1}) \cdots \sub^n(S^x_{m+n}) \cdots.$$
Hence for $m=-n$,  recalling that $\omega^x \in [W]$,  we obtain (looking at the right hand part)
$$\omega^{\frak{S}^{-n }x} _+=  c^x_{-n} S^x_{-n} \sub(S^x_{-n+1}) \cdots \sub^{n-1}(S^x_{-1}) \sub^n(W_1\cdots W_{\dd}) \cdots $$
By construction, $U^n$ is a prefix of this word (shifted once). It remains to compute its length. We observe that $|U^0| = \tau_0 -1$. But  $|U^{n+1}| = |S^x_{-(n+1)}| + \dd |U^n| + \tau_{-(n+1)} + m -1 = \dd |U^n| + \dd +\tau_{-(n+1)}$. Hence if we assume $|U^n| = \dd^n (\sum_{k=-n}^{0} \tau_k d^k )-1 $, we obtain $|U^{n+1}| =  \dd^{n+1} (\sum_{k=-n}^{0} \tau_k d^k ) - \dd + \dd +\tau_{-(n+1)} = d^{n+1} \sum_{k=-(n+1)}^{0} \tau_k \dd^k -1 (=\lfloor \dd^{n+1} t \rfloor -1).$
In view of the definitions of $t_+^n$ and $t_-^n$, combining this discussion with equations  \eqref{equation:ergo-sum1} and \eqref{equation:ergo-sum2}, we get the statement of the lemma.

 \end{proof}

To pass from limit theorems for the Markov chain to those for the ergodic integral, we need to combine Lemma \ref{lem:link-sums-automaton} with the following simple observation.
\begin{prop}
\label{varcompar}
Let $X_1^{(n)}$ and $X_2^{(n)}$, $n\in {\mathbb N}$,
be two sequences of random variables on a probability space $(\Omega, {\mathbb P})$
satisfying the conditions:
\begin{enumerate}
\item there exists a constant $K$ such that for all $n$ we have the inequality $|X_1^{(n)}-X_2^{(n)}|\leq K$ almost surely;
\item $Var( X_2^{(n)})\to\infty$ as $n\to\infty$.
\end{enumerate}
Then
$$
\lim\limits_{n\to\infty}\frac{Var( X_2^{(n)})}{Var( X_1^{(n)})}=1.
$$
\end{prop}
\begin{proof}
First, note that we only need to consider the case when ${\mathbb E}X_1^{(n)}={\mathbb E}X_2^{(n)}=0$, since if $X_1^{(n)}$ and $X_2^{(n)}$ satisfy the assumptions of the Proposition, then so do $X_1^{(n)}-{\mathbb E}X_1^{(n)}$ and $X_2^{(n)}-{\mathbb E}X_2^{(n)}$. If ${\mathbb E}X_1^{(n)}={\mathbb E}X_2^{(n)}=0$, then we have
$$
|{\mathbb E}(X_1^{(n)})^2-{\mathbb E}(X_2^{(n)})^2|\leq K\left(K+{\mathbb E}|X_1^{(n)}|\right)\leq K\left(K+\sqrt{Var( X_2^{(n)})}\right)
$$
by the Cauchy-Buniakovsky-Schwarz inequality, and the proof is complete.
\end{proof}

\section{Markov approximation and the proof of Theorem  \ref{fluctuations-periodic} }
\label{proof-fluctuations-periodic}
The following construction plays a basic role in the sequel. 
\subsection{Markov chains}  Let $t \in \RR_+^*$. Let $n_0$ be the unique integer such that $\dd^{n_0}t \in [1;\dd)$ and set $t=\sum_{k \leq n_0} \tau_k \dd^k$. For now, we work with $n_0 = 1$. Let $(\Omega, \F, P)$ be a probability space. We construct a random variable $I$ on $\Omega$ valued in $\A^\dd$ with distribution
$$P(I = W)  = \nu(\left\{x \in Y : \omega^x_0 \cdots \omega^x_{\dd} = W \right\}).$$
We define a map $g^t_0 : \A^\dd \to \RR$ by $g_0(W) = f(W_1 \cdots W_{k_t-1})$. Let now $(m_n)_{n \geq 1}$ be a sequence of iid random variables uniformly distributed   in $\{1, \ldots, \dd\}$ independent of $I$. We define a process valued in $\E = \A \times \A^2$ by setting $X_0 = (I_1, I_k I_{k+1})$ and recursively $X_{n+1}=(a^{n+1},V^{n+1})$ to be the end  vertex of the edge labelled $m_n$ of the automaton $A_{\tau_{-n}}$ starting from $X_n=(a^n,V^n)$.
Finally, we set

\begin{equation} \label{gn}
g^t_n(X_n,m_n) = f(S^{a^n,m_n}) + f(P^{V^n,m_n+\tau_{-n}}).
\end{equation}

\begin{lem}
The sequence
$(X_n)_{n \geq 0}$ is a Markov chain.

We have the following commutative diagram

$$\begin{array}[c]{ccc}
\E^\NN&\stackrel{\frak{S_{t}}}{\rightarrow}&\E^\NN\\
\downarrow\scriptstyle{\pi}&&\downarrow\scriptstyle{\pi}\\
Y&\stackrel{\frak{S}}{\rightarrow}&Y
\end{array}$$
where $\frak{S_{t}}$ is the shift on the Markov chain. The  probability measure $P$ invariant by the Markov chain projects to $\nu$.

\end{lem}

\begin{proof}
$(X_n)_{n \geq 0}$ is a Markov chain since $(a_{n+1}, V_{n+1})$ only depends on $(a_{n}, V_{n})$.
Moreover, there is a projection
$\pi$ between infinite paths of $(X_n)$ and $Y$. This map consists in forgetting the second variable $V^n$. Formally, the map is defined by : $m_n = |P_n^x| +1$.
 To check that the image measure is $\nu$ under the projection,  it suffices to observe that $\sum_{V \in \A^2} P(X_0 = (a,V)) = \nu(\{x \in Y :  \omega^x_0 = a\}) $ which is obvious.
\end{proof}

We claim that
\begin{lem} \label{lem:estimate-markovchain}
$$f(U^n) = g_0(X_0) + \sum_{k=1}^{n} g_k(X_k,m_k). $$
\end{lem}
\begin{proof}
It is a straightforward computation. Check that $f(U^0) = g_0(X_0)$. Then write recursively $U^{n+1}=S_{-(n+1)}^x \sub(U^n) V^n_1 \cdots V^n_{\tau_{-n}+m-1}$. Hence $f(U^{n+1}) - f(U^n) = f(S_{-(n+1)}^x) + f(V^n_1 \cdots V^n_{\tau_{-n}+m-1})$.
\end{proof}

Now, we assume that $t$ has a periodic expansion thus the Markov chain is homogeneous. To simplify notations,  we will forget the dependence in $t$. Up to a change of basis, we assume that the expansion of $t$ has period 1, all the digits are equal, i.e. for all $k$, $\tau_{k} = \tau$.

%

\subsection{Proof of Theorem \ref{fluctuations-periodic}}

$\mbox{ }$

The Markov chain $(X_{n})$ is endowed with  an initial probability measure $\mu_{0}$ defined on
$\A  \times \A^2$ and a transition matrix $Q$.

\begin{lem}
Let $B = \{ V \in \A^2, \nu([B]) >0\}$.
The Markov chain $(V_{n})$ is recurrent on $\A \times B$.
\end{lem}

\begin{proof}
The proof is straightforward since the substitution is primitive.
\end{proof}

Observe that the processes $(a_n)$ and $(V_n)$ are also Markov chains. Let $\tilde{\mu}$ be an invariant measure for the Markov chain $(X_{n})$.

\begin{lem} \label{invariant-probas}
The initial measure $\tilde{\mu_{0}}$ of the Markov chain $(X_{n})$ projects on the first coordinates on the invariant measure for the Markov chain $(a_{n})$ and on the last coordinate on the invariant measure for the Markov chain $(V_{n})$.
\end{lem}

\begin{proof}
The initial measure satisfies:
$$\mu_{0}(a,m,V) = P(I_{1} = a, I_{k_{t}I_{k_{t}+1}} = V) P(m_{0} = m).$$ It projects on $\mu_{0}^1 (a) = \nu([a])$ and on
$\mu_{0}^2 (V) = \nu([V])$. These measures are obviously invariant under the Markov chains.
\end{proof}

An important lemma is the following:

\begin{lem} For every integer $n$ and for every probability measure $\tilde{\mu}$  invariant by the chain $(X_{n})$, we have
$E[g_{n}(X_{n})] = 0$ and $E_{\tilde{\mu}}[g_{0}(X_{0})] = 0$.
\end{lem}

\begin{proof}

We recall that $E_{\nu}[\int_{0}^{\dd^nt} f \circ h_s(x) ds] = 0$ since $E[f] = 0$ by hypothesis.
This does not {\it a priori} imply the conclusion of the Lemma since
 the Markov chain $(X_{n})$ is not always recurrent.
 We must exclude the possibility the averages
 on different ergodic components compensate each other.

By Lemma  \ref{lem:estimate-markovchain}, there exists a constant $K$ such that, independently on $n$, 
\begin{equation} \label{upper bound}
E[\sum_{k=1}^{n} g_k(X_k,m_k)] <K. 
\end{equation}
Since expectation is a linear operator, and in view of the definition of $g_{n}$ (recall equation \eqref{gn}), we can decompose
\begin{equation} \label{expectation-linear}
E[g_{n}(X_{n}, m) ] = E[g_{n}^1(a_{n}, m) ] + E[g_{n}^2(V_{n}, m) ].
\end{equation}

Let $\tilde{\mu}$ be any invariant measure under $(X_{n})$. We have:
$$E_{\tilde{\mu}} [g_{0}(X_{0},m)] = E_{\tilde{\mu}} [g_{n}(X_{n},m)] = E_{\tilde{\mu}} [g_{n}^1(a_{n}, m) ] + E_{\tilde{\mu}} [g_{n}^2(V_{n}, m) ].$$
By Lemma \ref{invariant-probas} and equation \eqref{expectation-linear} this yields 
$$E_{\tilde{\mu}} [g_{0}(X_{0},m)] = E[g_{n}^1(a_{n}, m) ] + E[g_{n}^2(V_{n}, m) ] = E[g_{n}(X_{n}, m) ] .$$

Therefore, $E[g_{n}(X_{n}, m) ] $ is independent on $n$. It follows from \eqref{upper bound}  that $E[g_{n}(X_{n}, m) ] = 0 $ and  $E_{\tilde{\mu}}[g_{0}(X_{0}] = 0$.

\end{proof}

We just proved that  for every recurrent class the mean of the  invariant measure supported by this class is equal to zero.
Thus the central limit theorem for stationary Markov chains implies the random variable
$\frac{1}{\sqrt{n}}\sum_{k=1}^n g_{k}(X_{k},m_{k})$ has a limiting distribution. The contribution of each recurrent class is a normal law.
\begin{remark}
Observe that, as a consequence,  for every positive $\varepsilon$,
$$\frac{1}{n^{\frac{1}{2}+\varepsilon}}\sum_{k=1}^n g_{k}(X_{k},m_{k})$$
converges in distribution to 0.
\end{remark}
To finish the proof of Theorem  \ref{fluctuations-periodic}, it is enough to show:
\begin{equation}
 \lim_{ n \to \infty} \frac{1}{\sqrt{n}} \int_{0}^{\dd^nt} f \circ h_s(x) ds  \stackrel{(d)}{=}
 \lim_{ n \to \infty}  \frac{1}{\sqrt{n}} \sum_{k=1}^{n} g_k(X_k,m_k).
\end{equation}
In other words :
\begin{prop}
For all bounded continous $\varphi : \RR \to \RR$,
\begin{equation*}
 \lim_{ n \to \infty} \int_Y \varphi \left( \frac{1}{\sqrt{n}} \int_{0}^{\dd^nt} f \circ h_s ds\right) \nu  =
\lim_{ n \to \infty}  E\left( \varphi\left( \frac{1}{\sqrt{n}} \sum_{k=1}^{n} g_k(X_k,m_k) \right)\right) .
\end{equation*}
\end{prop}
\begin{proof}
By Lemmas \ref{lem:link-sums-automaton} and
\ref{lem:estimate-markovchain}, the quantities $ \int_{0}^{\dd^nt} f \circ h_s(\frak{S}^{-n} x) ds$  and $\sum_{k=1}^{n} g_k(X_k,m_k)$ are equal up to a  bounded term. As we divide by $\sqrt{n}$ a bounded term is negligible. Thus the limit distribution of  the sequences
$$\frac{1}{\sqrt{n}} \int_{0}^{\dd^nt} f \circ h_s(\frak{S}^{-n}  x) ds \, \hbox{ and } \frac{1}{\sqrt{n}} \sum_{k=1}^{n} g_k(X_k,m_k)$$ are the same.

By Lemma \ref{invariance-under-frakS}, the law of
$$ \int_{0}^{\dd^nt} f \circ h_s( x) ds$$ is equal to the law of  $$ \int_{0}^{\dd^nt} f \circ h_s(\frak{S}^{-n} x) ds.$$ This proves the proposition.  \end{proof}
Together with Proposition \ref{varcompar}, this implies Theorem  \ref{fluctuations-periodic}.


\section{Proof of Theorem \ref{fluctuations-generic}} \label{proof-thm-generic}


\subsection{Assumptions on the Markov chain}

In the beginning of Section \ref{section:automata}, to a substitution $\sigma$ we  assigned a family of automata
$\{ A_{i},  i \in \{1, \dots, \}$ in such a way that to any $t \in [0,1]$ with decomposition in basis $d$ (set $\A=\{0, \ldots, d-1\}$) given by $t=\sum_{n >  0} \tau_n d^{-n}$  we associate a sequence of automata in the family, $A_{\tau_{0}}, A_{\tau_{1}}, \ldots$. To such $t$ we also associate a sequence of functions $(g^t_k)$ defined on the state space $E$ of the underlying Markov chain ($g^t_k = g_{\tau_k}$). When $t$ has a periodic expansion in basis $d$, the sequence of automata is periodic and hence the associated Markov chain is a standard homogeneous Markov chain.

Here we assume  that there is such a $t_0$ for which the Markov chain is {\em aperiodic and has positive variance}.
Let us first explain the meaning of this assumption. We stress that it is not only an assumption on the Markov chain but also on the functions $(g^t_k)$.
Let  $\chi = \chi_1 \cdots \chi_p$  denote a word in $\A^p$ such that  $t_0 = \sum_{n < 0} (\sum_{i=1}^p \chi_i d^{-i}) d^{-np}$.
First of all, we assume that the matrix of the Markov chain associated with $t_0$ is primitive  (without loss of generality, we can assume that the matrix has strictly positive entries). The positive variance asumption means that the sequence of functions $g^{t_0}$ does not yield a coboundary. Hence there are two cycles on which it takes two different values.

In slightly different words,  it implies that there are two paths $u$ and $v$ in $E^{p+1}$  (with $u_0=v_0$ and $u_p=v_p$) and some $\epsilon>0$
such that
\begin{equation}
\label{deuxchemins}
\left| \sum_{k=1}^{p} (g_{\chi_{k}} (u_{k})) - \sum_{k=1}^{p} (g_{\chi_{k}}(v_{k})) \right|  \geq \epsilon > 0,
\end{equation}
We claim that under this assumption, for Lebesgue almost all $t$,
\begin{equation} \label{intermediate-convergence1}
\frac{\sum_{k=0}^{n-1} (g^t_k(X_k) - E[ g^t_k(X_k)] )}{\sqrt{V( \sum_{k=0}^{n-1} g^t_k(X_k)))} }\Rightarrow N(0,1).
\end{equation}
As we shall see more precisely in Section \ref{proof-fluctuations-generic}, by Lemma \ref{lem:link-sums-automaton} and Lemma \ref{lem:estimate-markovchain}, convergence \eqref{intermediate-convergence1} implies Theorem \ref{fluctuations-generic}.
\subsection{Choice of a full-measure set of admissible sequences}
Observe that the Lebesgue measure on $[0,1]$ maps onto the uniform product measure on $\A^\NN$ by the map sending a real number to its expansion in basis $d$. Hence, with probability $1$, the expansion of $t$ contains infinitely many occurrences of $\kappa=\chi^3$ (we stress that we ask for three successive occurrences of $\chi$).

Denote $(q_n)$ the sequence of occurrences of $\kappa$ in the expansion of $t$ and $m_n=q_{n}+3p$ (recall $p=|\chi| = |\kappa|/3$ ; fix $m_0=0$). The sequence $(m_n-m_{n-1})$ is a sequence of independent identically distributed random variables with geometric law.  We claim that it does not grow too fast.

 We use Borel-Cantelli Lemma to say that the inequality $B_n = \{m_{n} - m_{n-1} > n^{\frac{1}{4}} \} $ almost surely  holds only for a finite number of $n \in \NN$ since  $\sum_{n \geq 1} \PP(B_n) = \sum_{n \geq 1}  c \, \theta^{n^{1/4}} < \infty$. Let $J \subset [0,1]$ be the set of probability $1$ on which  the expansion of $t$ contains infinitely many occurrences of $\kappa$ and $B_n$ is true a finite number of times. Observe that for all $t \in J$, there is $C$ (depending on $t$)  such that for all $n \in \NN$,
 \begin{equation}
 \label{emen}
 m_{n} - m_{n-1} <  C n^{\frac{1}{4}}.
\end{equation}
\vip


\subsection{A Central Limit Theorem for non homogeneous Markov chains}
\subsubsection{The Dobrushin Theorem}
We recall
Dobrushin's central limit theorem for non-homogeneous Markov chains (see \cite{Dob}). We follow the exposition
by Sethuraman and Varhadan \cite{SV}. We only state the theorem for  Markov chains defined on a finite state space,
since that result is sufficient for our purposes.

For each $n\geq1$, let $\{X_i^{(n)} : 1\leq i \leq n \}$ be $n$ observations of a non-homogeneous Markov chain on a finite state space $E$ with transition matrices $\{\pi^{(n)}_{i,i+1} : 1 \leq i \leq n-1\}$ and initial distribution $\mu^{(n)}$.

Recall that, for discrete Markov chains,  the Dobrushin's ergodic coefficient of a transition matrix $\pi$ is defined as $\alpha(\pi) = 1 - \delta(\pi)$ where
$$\delta(\pi) = \sup_{a,b \in E ; c \in E} |\pi(a,c) -\pi(b,c)|. $$
Recall the well-known inequality:
\begin{equation} \label{eq:ergodic-coef}
\alpha(\pi) \geq \min_{a,c \in E} \pi(a,c).
\end{equation}

Moreover if $\eta$ is a continuous function its oscillation is  $Osc(\eta) = \sup_{\omega,\omega'} |\eta(\omega) - \eta(\omega')|$. If $\eta$ is a function defined on the states of a Markov chain, we have the property :  $Osc(\pi \eta) \leq \delta(\pi) Osc(\eta)$.

Let $\alpha_n = \min_{1\leq i \leq n-1} \alpha(\pi^{(n)}_{i,i+1})$. Let
$\{f_i^{(n)} : 1 \leq i \leq n\}$ be real valued functions on $E$ such that there exists some finite constants $C_n$ with $\sup_{1\leq i \leq n} \sup_{x \in X} {|f_i^{(n)}(X_i^{(n)})|} \leq C_n. $
Define, for $n\geq 1$, the sum
$$S_n = \sum_{i=1}^{n} f_i^{(n)} (X_i^{(n)}).$$
\begin{theo}[Theorem 1.1 in \cite{SV}]
If
$$\lim_{n \to \infty} C_n^2 \alpha_n^{-3} \left[\sum_{i=1}^{n} V(f_i^{(n)}(X_i^{(n)}))\right]^{-1} = 0,$$
then, we have the standard normal convergence
$$\frac{S_n - E[S_n]}{\sqrt{V(S_n)}} \Rightarrow N(0,1).$$
\end{theo}

\subsubsection{Decomposition for mixing}

Let now $t \in J$ be fixed. Most quantities depend on $t$ but, from now on,  we omit the sub/superscript.  Let $(\tau_n)_{n\geq 0}$ be the expansion of $t$. We denote $p_{n}$ the transition matrix associated with the automaton $A_{\tau_{n}}$ and $p_{n,m} = p_n \cdots p_{m-1}$ the transition matrix between  times $n$ and $m$.
For all $n$ consider the word $W_n = \tau_{m_{n-1}} \cdots \tau_{m_{n}-1}$.
This yields a decomposition of our sequence in words ending with $\kappa=\chi \chi \chi$. We define the (non homogeneous)  Markov chain  $Y_n = X_{m_{n-1}}$. We denote by $\pi_{n}$ its transition matrices, $\pi_{n}=p_{m_{n-1}, m_{n}}$. The function we want to control is a function of $(X_{m_{n-1}}, \ldots, X_{m_{n}-1})$.
Let us set $\tilde Y_n = (X_{m_{n-1}+1}, \ldots, X_{m_{n}-1})$. Observe that, conditionally to $Y_n$ and $Y_{n+1}$,  $\tilde Y_n$ is independent of $\{Y_i, i \in \NN \setminus \{n, n+1\}\}$ (and of $\{\tilde Y_i, i \in \NN \setminus \{n\} \}$). We introduce
$$f_n(Y_n, \tilde Y_n) = \sum_{k=m_{n-1}}^{m_{n}-1} g_k(X_k) - \EE [ \sum_{k=m_{n-1}}^{m_{n}-1} g_k(X_k)] .  $$

\begin{remark}
Ideally we would like to use Dobrushin's result about sums of observables of a non stationary Markov chain. The difficulty is that the functions do not depend only on the Markov chain itself. We could represent this quantity as a function of $Y_n$, $Y_{n+1}$ and another random variable independent of the chain. Still it is not enough to apply the theorem in its standard versions. The usual way to avoid this problem would be to build a new chain joining coordinates by pairs $((Y_n,Y_{n+1}))$; but this construction breaks down the lower bound for the ergodic coefficient. This does not change the result fundamentally but we have to explain how the proof has to be adapted.
\end{remark}

We state the main results of this section as follows.
We denote $F_n = f_n(Y_n, \tilde Y_n) $ and $S_n = \sum_{k=1}^{n} F_k$. The positive variance assumption allows to show that, for all $t \in J$,
\begin{lem}
\label{lem:variance}
There is $\gamma >0$ such that
\begin{equation}
\label{thevariance}
V\left(S_n \right) \geq \gamma n.
\end{equation}
\end{lem}
The proof is given in Section \ref{subsection-variance}. Relying on this result, we show that
 \begin{prop}
\label{clt1}
For all $t \in J$,
$$\frac{S_n}{\sqrt{V(S_n)}} \Rightarrow N(0,1)$$
\end{prop}

The remainder of this section is devoted to the proof of Proposition \ref{clt1}.   We are going to adapt Dobrushin's
Central Limit Theorem for inhomogeneous Markov chain (following the exposition by Sethuraman and Varhadan \cite{SV}) to our special case. But we first show how it yields Theorem \ref{fluctuations-generic}.


\subsection{Proof of Theorem \ref{fluctuations-generic}}
\label{proof-fluctuations-generic}
We are now in position to prove Theorem \ref{fluctuations-generic}. Let us fix $t \in J$.

\vip

Lemma \ref{lem:link-sums-automaton} and Lemma \ref{lem:estimate-markovchain} show that we can apply  Proposition \ref{varcompar}, to get equivalence of the variances $V_n$  introduced in the statement of Theorem  \ref{fluctuations-generic} and  $V( \sum_{k=0}^{n-1} g^t_k(X_k)))$. Since the $g_k^t$ belong to a  finite family of bounded functions, we easily get an upper bound for the variances: $V( \sum_{k=0}^{n-1} g^t_k(X_k)) \leq \sup_{k \geq 0} ||g^t_k||_\infty n\leq \alpha_2 n$ . For the lower bound, we have to be more careful. On the one hand, Lemma \ref{lem:variance} shows that
there is $\gamma >0$ such that $V\left(S_n \right) \geq \gamma n$.
On the other hand, we use the relationship (\ref{emen}) between $n$ and $m_n$ to observe that $m_n \leq n . C.n^{1/4} \leq C. n^{5/4}$. Hence $V( \sum_{k=0}^{m_n-1} g^t_k(X_k))) \geq \gamma (m_n/C)^{4/5}$ and (for $m_n \leq m < m_{n+1})$, $V( \sum_{k=0}^{m-1} g^t_k(X_k))) \geq \gamma (m_n/C)^{4/5} \geq \gamma ((m-m^{1/4})/C)^{4/5} \geq \alpha_1 n^{4/5}$. This proves statement $(i)$. 

\vip

We use again (\ref{emen}) and boundedness of the family $(g_k^t)$ to check that  for $m_n \leq m < m_{n+1}$, $|\sum_{k=0}^{m-1} g^t_k(X_k))) - S_n| \leq \hbox{const. } n^{1/4}$.
In view of the variance lower bound, Proposition \ref{clt1}  yields
\begin{equation} \label{intermediate-convergence}
\frac{\sum_{k=0}^{n-1} (g^t_k(X_k) - E[ g^t_k(X_k)] )}{\sqrt{V( \sum_{k=0}^{n-1} g^t_k(X_k)))} }\Rightarrow N(0,1).
\end{equation}
Lemma \ref{lem:link-sums-automaton} and Lemma \ref{lem:estimate-markovchain} show that (\ref{intermediate-convergence}) yields  statement (ii) of Theorem \ref{fluctuations-generic}.

\vip

We conclude recalling that $J$ has  full Lebesgue measure.  $\Box$


\subsection{Standard assumptions}  Even though our situation is very similar, we cannot apply directly Dobrushin's Theorem  because our sum runs on observables that are not functions of the state of the (mixing) Markov chain. Nonetheless, our argument follows the general outline of \cite{SV}. We  use the same interplay between boundedness, rate of mixing and variance. More specifically, in our case, the ergodic coefficient is constant (bounded away from $0$ ; so we need just to ensure that $C_n^2/V(S_n) \to 0$).

\vip

\subsubsection{Upper bound} The upper bound
$$C_n = \sup_{i\leq n} \sup_{z,\tilde z }{f_i(z,\tilde z)} \leq \sup_{i\leq n} (||g_i||_\infty  (m_{i+1} - m_i)) \leq C_1 .  n^{1/4}.$$
In other words,
$$C_n =  \sup_{i\leq n} F_i \leq C_1 n^{1/4}.$$

\vip

\subsubsection{Ergodic coefficient} The ergodicity coefficient is bounded from below uniformly because of the presence of the word $\kappa$ at the end which guarantees some mixing. Let us be more specific.

By inequality \eqref{eq:ergodic-coef} since all the words $W_n$ end with the word $\kappa$ (associated to a (power of a) primitive matrix),  $p_{q_n ,m_{n}}$ is uniformly bounded from below (by a constant $\alpha$), and we have
$$\pi_n(a,b) \geq \sum_{c \in A}  p_{m_{n-1},q_n}(a,c) p_{q_n,m_{n}}(c,b) \geq \sum_{c \in A} p_{m_{n-1},q_n}(a,c)  \,  \alpha  \geq \alpha. $$
From which we deduce that $\alpha(\pi_n) \geq \alpha$ uniformly in $n$. This holds for the Markov chain $(Y_n)_{n\geq 0}$.

\begin{remark}
\label{detail}
We stress a technical difference with \cite{SV}. The point is that, as in \cite{SV}, we have  $Osc(E[f(Y_j,\tilde Y_j)|Y_i]) \leq \alpha^{j-i} Osc(f(Y_j))$ ; but that we will need a bound on
$Osc(E[f_j(Y_j,\tilde Y_j)|Y_i,\tilde Y_i])$ where $\tilde Y_i$ depends on what happens between $m_{i-1}$ and $m_{i}-1$. Since we need only asymptotic decay, i.e. for $j-i$ large, we will simply shift by $1$ using that
$$Osc(E[f_j(Y_j,\tilde Y_j)|Y_i,\tilde Y_i])
\leq Osc(E[f_j(Y_j,\tilde Y_j)| Y_{i+1}] ).$$
\end{remark}
\vip

\subsubsection{Variance} The variance of each elementary contribution is uniformly bounded from below. Furthermore the sum of the variances of the first $n$ terms is increasing rapidly enough with $n$.  In the usual setting the last statement follows from the first one using assumption 2. In our case, it seems more convenient to prove directly the second statement. That is the object of Lemma \ref{lem:variance} which claims that there is $\gamma >0$ such that $V\left(S_n \right) \geq \gamma n$ and whose proof is postponned to  Section \ref{subsection-variance}.

\subsubsection{Outline of the argument}
The crucial point for the proof (following \cite{SV}) is the interplay between these three quantities. Roughly speaking the variance must grow fast enough to kill unboundedness and lack of ergodicity, as expressed by equation $(1.3)$ in \cite{SV}. In our case the ergodic coefficient is independent on $n$ so it is enough to have the variance growing faster than the upper bound. To draw the parallel let us write :
$$ \lim_{n \to \infty} C^2_n \alpha_n^{-3} V_n^{-1} =  \lim_{n \to \infty} C^2 n^{1/2} \rho^{-3}   \frac{1}{\gamma n} = 0.$$
But we stress again that the situation is slightly more tricky because the observable we are interested in is not a simple function of the states of the Markov chain for which the assumptions are fullfilled ($f_n$ depends on $Y_n$, but also on $\tilde Y_n$). That is why we treat the variance separately. For the remainder of the proof, we can not just quote \cite{SV} but it appears that we can follow exactly the same lines. To construct the Martingale approximant, we take $\tilde Y_n$ into account. To prove that the assumptions of the main Martingale differences CLT arguments are fulfilled, we use the same ideas based on the control of the oscillations of the conditionnal expectations of $f_n(Y_n,\tilde Y_n)$ which rely essentially on the summability of the long range correlations and on the sublinearity of the upper bound (together with the linear growth of the variance).

\subsection{Proof of Proposition \ref{clt1}}
The central idea is to use the following standard CLT (implied by Corollary 3.1 in \cite{HH}) for martingale differences :  if $M^{(n)}_k = \sum_{l=2}^{k} \xi^{(n)}_l$ is a martingale  with respect to a filtration $\F_k$ and if,
\begin{equation}
\label{negli}
\max_{1 \leq i \leq n} ||\xi^{(n)}_i||_\infty \to 0, \, \, \hbox{ and, }
\end{equation}
\begin{equation}
\label{carreint}
\sum_{i=1}^{n} E[(\xi^{(n)}_i)^2]|\F_{i-1}] \to 1 \hbox{ in } L^2,
\end{equation}
then,
$$M_n^{(n)} \Rightarrow N(0,1).$$

The quantities whose fluctuations we study are not a martingale mainly because increments  are not independent for a Markov chain in that they depend on the present state. It is rather  standard to bypass this difficulty by introducing
\begin{equation}
\label{defZ}
 Z^{(n)}_k = \sum_{i=k}^{n} E[ f_i(Y_i, \tilde Y_i) \, | \, Y_k, \tilde Y_k].
 \end{equation}
Then considering the scaled differences
$$\xi_k^{(n)} = \frac{1}{\sqrt{V(S_n)}} \left(  Z^{(n)}_k - E[Z^{(n)}_k \, | \, Y_{k-1}, \tilde Y_{k-1}] \right), $$
we define the (family of) processes  $M^{(n)}_k = \sum_{i=2}^{k} \xi^{(n)}_i$ and observe that the processes $(M^{(n)}_k)_{1\leq k \leq n}$  are martingales  with respect to the filtrations $\{\F_k, 1 \leq k \leq n \}$, where $\F_k = \sigma\{ Y_l , \tilde Y_l : 1 \leq l \leq k \}$.
Indeed, since $ \xi^{(n)}_{k}$ is $\F_k$-measurable, for $k\leq n-1$,
\begin{eqnarray*}
E[M_{k+1}^{(n)} | \F_{k}] 
&=& \sum_{i=2}^{k} \xi^{(n)}_i   + E [ \xi^{(n)}_{k+1} |\F_{k}] \\
&=&M_k^{(n)} + (\sqrt{V(S_n)})^{-1} (E\left[ Z^{(n)}_{k+1} - E[Z^{(n)}_{k+1} \, | \, Y_{k}, \tilde Y_{k}]   \,  | \, \F_{k}\right] ) \\
&=&M_k^{(n)} + (\sqrt{V(S_n)})^{-1} (E[Z^{(n)}_{k+1}| \F_k]  - E[Z^{(n)}_{k+1} \, | \, Y_{k}, \tilde Y_{k}]    |\F_{k}] ) \\
&=&M_k^{(n)} .
 \end{eqnarray*}
Observing that
\begin{eqnarray*}
E[Z^{(n)}_{k+1} \, | \, Y_k, \tilde Y_k] &= &E[      \sum_{i=k+1}^{n} E[ f_i(Y_i, \tilde Y_i) \, | \, Y_{k+1}, \tilde Y_{k+1}]     \, | \, Y_k, \tilde Y_k]\\
&= &    \sum_{i=k+1}^{n}  E[     E[ f_i(Y_i, \tilde Y_i) \, | \, Y_{k+1}, \tilde Y_{k+1}]     \, | \, Y_k, \tilde Y_k]\\
&= &    \sum_{i=k+1}^{n}   E[ f_i(Y_i, \tilde Y_i)   \, | \, Y_k, \tilde Y_k], \\
\end{eqnarray*}
we can put (\ref{defZ})  the other way round to obtain,
for $k\leq n-1$,
$$f_k(Y_k,\tilde Y_k) = Z^{(n)}_k - E[Z^{(n)}_{k+1} \, | \, Y_k, \tilde Y_k],$$
for $k=n$, $f_n(Y_n,\tilde Y_n) = Z^{(n)}_n.$ We obtain, for the sums we are interested in:
$$ S_n = \sum_{k=2}^{n} \left( Z^{(n)}_k - E[Z^{(n)}_k \, | \, Y_{k-1}, \tilde Y_{k-1}] \right) + Z_1^{(n)},  $$
and for their variances (using the orthogonality of the martingale increments): 
\begin{equation}
\label{detailvariance}
V(S_n) = \sum_{k=2}^{n} V\left( Z^{(n)}_k - E[Z^{(n)}_k \, | \, Y_{k-1}, \tilde Y_{k-1}] \right)  + V(Z_1^{(n)}. 
\end{equation}
Hence, we can approximate $S_n / \sqrt{V(S_n)}$ by $M^{(n)}_n$ (see Step 5.) and use the above argument to conclude, provided we  check both assumptions (\ref{negli}) and (\ref{carreint}) (as well as approximation).

\vip
To do so, we will follow the proof of \cite{SV}.

\vip

{\bf Step 1. } We prove inequalities  (playing the role of Lemma 3.1 in \cite{SV}).  We will make repeated use of the property :  $Osc(\pi \eta) \leq \delta(\pi) Osc(\eta)$.

Let $1\leq  i  < j \leq n$. As $|| F_j^{(n)}||_\infty \leq C_n$, its oscillation $Osc(F_j) \leq 2 C_n$.
From the control we have on the mixing coefficient (see Remark \ref{detail}),
$$Osc( E[F_j^{(n)} | \F_i]) \leq Osc( F_j^{(n)}) \delta(\pi_{i+1,j}) \leq 2 C_n (1-\alpha)^{j-i-1}.$$
For centered random variables, $||||_\infty$ is bounded by $2 Osc()$, so that, for all $1 \leq i<  j \leq  n$,
\begin{equation}
\label{ineq1}
||E[F_j|\F_i ||_\infty \leq 2 Osc( E[F_j^{(n)} | \F_i]) \leq \hbox{const. } n^{1/4} (1-\alpha)^{j-i}.
\end{equation}


The same arguments also yield the similar
\begin{equation}
\label{ineq2}
Osc(E[F_j^2|\F_i]) \leq  \hbox{const. } n^{1/2} (1-\alpha)^{j-i}.
\end{equation}

{\bf Step 2. } Putting together inequality (\ref{ineq1})  and the variance lower bound Lemma~\ref{lem:variance} (corresponding to Proposition 3.2 \cite{SV}), we obtain an analog of  Lemma 3.2 \cite{SV}.

Using inequality (\ref{ineq1}) and $EF_i =0$, we see that
$$||Z^{(n)}_k ||_\infty \leq \sum_{i=k}^n || E[F_i^{(n)} | \F_k ]  ||_\infty \leq \hbox{const. } n^{1/4}  \sum_{i=k}^n  (1-\alpha)^{i-k} \leq 2 C_n/\alpha.$$
Then by the variance lower bound (Lemma~\ref{lem:variance}),
\begin{equation}
\label{norminf}
\sup_{1\leq k \leq n}  \frac{||Z^{(n)}_k||_\infty}{\sqrt{V(S_n)}}  \leq  \hbox{const. }\frac{n^{1/4}}{\alpha \sqrt{\gamma n}},
\end{equation}
which is, under our assumptions $o(1)$.

\vip

This will yield the assymptotic equivalence of $S_n/\sqrt{V(S_n)}$ and $M^{(n)}_n$ as well as negligibility (\ref{negli}) of the differences $\xi^{(n)}_k$.

\vip

{\bf Step 3. } Now let us set $H_i = E[(\xi^{(n)}_i)^2]|\F_{i-1}] $. Before all, we observe that
$$E[   \sum_{i=1}^{n} H_i ] = \sum_{i=1}^{n}E [ E[(\xi^{(n)}_i)^2]|\F_{i-1}] ] =   \sum_{i=1}^{n} E[(\xi^{(n)}_i)^2] =1 + o(1)$$
in view of (\ref{detailvariance}) and (\ref{norminf}).
In the next two steps, we are going to prove that $ \sum_{i=1}^{n} H_i  \to 1 \hbox{ in } L^2$. This first step is rather general and correspond to Lemma 3.3 in \cite{SV}. Observe that, since  $E[   \sum_{i=1}^{n} H_i ]=1 +o(1)$, the $L^2$-norm
$$ E[ ( \sum_{i=1}^{n} H_i  - 1 )^2] = E[ ( \sum_{i=1}^{n} H_i)^2] -2 E[   \sum_{i=1}^{n} H_i ] +1  =   E[ ( \sum_{i=1}^{n} H_i)^2] -1 + o(1)$$
tends to $0$ if and only if  $ E[  \left( \sum_{i=1}^{n} H_i \right)^2]  \to 1$.
Following the proof of Lemma 3.3 in \cite{SV} we write
$$ E[ ( \sum_{i=1}^{n} H_i)^2] =   \sum_{i=1}^{n} E[H_i^2] + 2 \sum_{i=1}^{n-1}E[H_i (\sum_{j=i+1}^nH_j)],  $$
and,
$$
 ( \sum_{i=1}^{n} E[H_i])^2 =  \sum_{i=1}^{n} E[H_i^2] + 2    \sum_{i=1}^{n-1} E[H_i] E[ (\sum_{j=i+1}^nH_j)].$$
Taking the difference and recalling $E[   \sum_{i=1}^{n} H_i ]=1 +o(1)$, we deduce that
$$ E[ ( \sum_{i=1}^{n} H_i  )^2] - 1 = 2 \sum_{i=1}^{n-1}E[H_i\left(  \sum_{j=i+1}^nH_j  -  E[ \sum_{j=i+1}^nH_j ]\right)] + o(1).$$

We observe that, since $H_i$ is measurable with respect to $\F_{i-1}$,
$$E[H_i \sum_{j=i+1}^nH_j ] = E[H_i E[\sum_{j=i+1}^nH_j| \F_{i-1}] ].$$
Hence,
\begin{eqnarray*}
| E[ ( \sum_{i=1}^{n} H_i  )^2] - 1|   & \leq & \sum_{i=1}^{n-1}E[H_i\left(  E[ \sum_{j=i+1}^nH_j | \F_{i-1}]  -  E[ \sum_{j=i+1}^nH_j ]\right)]  +o(1)\\
& \leq &\sup_{1 \leq i \leq n-1} Osc(E[ (\sum_{j=i+1}^nH_j)| \F_{i-1}])  \left( \sum_{i=1}^{n-1}E[H_i] \right) +o(1).
\end{eqnarray*}

\vip

{\bf Step 4. } Hence it remains only to prove
\begin{equation}
\label{oscillations}
\sup_{1 \leq i \leq n-1} Osc(E[ (\sum_{j=i+1}^nH_j)| \F_{i-1}]) = o(1). 
\end{equation}
 First we write
 \begin{eqnarray*}
E[ (\sum_{j=i+1}^nH_j)| \F_{i-1}]) &=&  E[\sum_{j=i+1}^n    E[(\xi^{(n)}_j)^2]|\F_{j-1}]  | \F_{i-1}])   \\
&=&  E[\sum_{j=i+1}^n   (\xi^{(n)}_j)^2 |\F_{i-1}] \\
\end{eqnarray*}
and use the  orthogonality property (martingale) of the increments ($E[\xi^{(n)}_r \xi^{(n)}_s | \F_u] = 0$ provided $r>s>u$) :
\begin{eqnarray*}
E[ (\sum_{j=i+1}^nH_j)| \F_{i-1}])
&=&  E[(\sum_{j=i+1}^n   \xi^{(n)}_j)^2 |\F_{i-1}]. 
\end{eqnarray*}

Then we write all this in terms of $Z$ and $F$. We recall that for $k \leq n-1$,  $F_k = Z^{(n)}_k - E[Z^{(n)}_{k+1}|\F_k]$, while $F_n = Z^{(n)}_n$.
\begin{eqnarray*}
\sum_{j=i+1}^n   \xi^{(n)}_j  &=& \sum_{j=i+1}^n Z^{(n)}_j - E[Z_{j}^{(n)} | \F_{j-1}] \\
&=& V(S_n)^{-1}  (  \sum_{j=i+1}^{n-1} F_j+ E[Z^{(n)}_{j+1}|\F_j]  - E[Z_{j}^{(n)}  | \F_{j-1}]) + F_n -  E[Z_{n}^{(n)}  | \F_{n-1}] ) \\
&=& V(S_n)^{-1} ( \sum_{j=i+1}^n F_j - E[Z_{i+1}^{(n)} | \F_i]). 
\end{eqnarray*}
Hence,
\begin{eqnarray*}
E[ (\sum_{j=i+1}^nH_j)| \F_{i-1}]) &=& V(S_n)^{-1} E[( \sum_{j=i+1}^n F_j - E[Z_{i+1}^{(n)} | \F_i])^2 | \F_{i-1}]  \\
&=&   V(S_n)^{-1} E[ (\sum_{j=i+1}^n F_j)^2 | \F_{i-1} ]     -  V(S_n)^{-1} E[ E[ Z_{i+1}^{(n)} |\F_i]^2 | \F_{i-1}]. 
\end{eqnarray*}

The last term is bounded by $\sup_{2 \leq i \leq n-1} V(S_n)^{-1} ||Z^{(n)}_{i+1}||_\infty^2 = o(1)$; hence its oscillation  is also uniformly $o(1)$.
For the first term, we write
$$Osc(V(S_n)^{-1}  E[ (\sum_{j=i+1}^n F_j)^2 | \F_{i-1} ] ) \leq V(S_n)^{-1}  \sum_{i+1 \leq j,m \leq n} Osc(E[F_j F_m | \F_{i-1}]).$$
But, in the same spirit as for the proof of inequality (\ref{ineq1}), we have, for all $1 \leq l < i<  j \leq  n$, 
\begin{eqnarray}
Osc(E\left[  F_i^{(n)}   E[F_j^{(n)}   | \F_i] | \F_l \right] ) &\leq& (1-\alpha)^{i-l-1}  Osc(F_i^{(n)}  E[ F_j^{(n)} | \F_i]  ) \nonumber \\
&\leq &   (1-\alpha)^{i-l-1}    \left( Osc(   F_i^{(n)})     || E[ F_j^{(n)} | \F_i]   ||_\infty   \right. \nonumber \\
&& \hspace{3cm} + \left. ||  F_i^{(n)}||_\infty  Osc(E[ F_j^{(n)} | \F_i] ) \right) \nonumber \\
\label{ineq3}
&\leq & \hbox{const. } n^{1/2} (1-\alpha)^{i-l} (1-\alpha)^{j-i} . \nonumber
\end{eqnarray}
\vip
It follows that, for all $i \leq n$,
$$Osc(E[ (\sum_{j=i+1}^nH_j)| \F_{i-1}]) \leq \hbox{const. } V(S_n)^{-1}  n^{1/2} /\alpha^2 = o(1).$$

\vip

\vip
{\bf Step 5. } Conclusion. On the one hand, Step 2. shows that almost surely,
$$\left| M^{(n)}_n - \frac{S_n}{\sqrt{V(S_n)}}\right| \leq \frac{||Z^{(n)}_1||_\infty}{\sqrt{V(S_n)}} = o(1). $$
On the other hand, Steps 1. and 2. show that  (\ref{negli}) 
holds while Steps 3. and 4. show that (\ref{carreint}) holds so that $M_n^{(n)} \Rightarrow N(0,1)$. This concludes the proof of Proposition \ref{clt1} $\Box$

\subsection{The variance} \label{subsection-variance}
We end the proof by giving a direct proof for Lemma \ref{lem:variance} about the lower bound of the variance.

\vip

{\em Step 1. } To start with, we claim that there is $\gamma>0$ such that, for all $i\geq 1$,
\begin{equation}
\label{var1}
E\left[ \left( f_i(Y_i, \tilde Y_i) - E\left[ f_i(Y_i, \tilde Y_i) | Y_i, Y_{i+1} \right] \right)^2  \, | \, Y_i, Y_{i+1} \right] > \gamma.
\end{equation}
To simplify we write $F_i = f_i(Y_i, \tilde Y_i) $ and  $G_i = F_i - E[F_i|Y_i,Y_{i+1}]$ so that (\ref{var1}) becomes :
$$E[ G_i ^2 \, | \, Y_i, Y_{i+1} ] > \epsilon.$$
To prove the claim, we fix $i \geq 1$. We set $m=m_{i-1}$, $M=m_i$, $p=|\chi|$,  $q=M_i-3p$, $r=M-2p$, $s=M-p$,  and $\ell = M-m$. We stress that with these notations,
$\tau_q \cdots \tau_{r-1} = \tau_r \cdots \tau_{s-1} = \tau_s \cdots \tau_{M-1} = \chi$ and recall that the matrix associated with $\chi$ is primitive.
\vip
For $w \in E^{p+1}$, we set $\E_w = \{X_r = w_0, \ldots , X_s = w_p\}$. Elementary considerations about Markov chains and primitivity of the matrices associated to $\chi$ yield existence of $\epsilon_1 >0$ and $\epsilon_2>0$ such that, for all $w \in E^{p+1}$,   $P(X_r = w_0, X_s = w_p \, | \, Y_i, Y_{i+1} ) > \epsilon_1$,
$P(\E_{w} \, | \, X_r=w_0, X_s=w_p) > \epsilon_2$.
\vip
We also observe that (\ref{deuxchemins}) shows that that there are two paths $u$ and $v$ in $E^{p+1}$  (with $u_0=v_0$ and $u_p=v_p$)
such that
$$\left| \sum_{k=r}^{s-1} (g_k(u_k)) - \sum_{k=r}^{s} (g_k(v_k)) \right|  \geq \epsilon > 0. $$

\vip

We observe that $E[ G_i ^2 \, | \, Y_i, Y_{i+1} ]  \geq E[ \indiq_{\E_u \cup \E_v} G_i ^2 \, | \, Y_i, Y_{i+1} ] $, so that, since $\E_u \cap \E_v = \emptyset$,
$E[ G_i ^2 \, | \, Y_i, Y_{i+1} ]  \geq E[ \indiq_{\E_u} G_i ^2 \, | \, Y_i, Y_{i+1} ]  +  E[ \indiq_{\E_v} G_i ^2 \, | \, Y_i, Y_{i+1} ]. $

\vip
By definition,
$$E[ G_i ^2 \, | \, Y_i, Y_{i+1} ] = \sum_{\tilde z \in A^{{\ell}}}  G_i(z_1\cdots z_\ell)^2 P(\tilde Y_i = \tilde z | Y_i, Y_{i+1}). $$
In the same language,
$$E[ \indiq_{\E_u} G_i ^2 \, | \, Y_i, Y_{i+1} ]  \geq \sum_{w \in E^{r}, z \in E^p} G_i(wuz)^2 P(X_m \cdots X_M = wuz |  Y_i, Y_{i+1}). $$
Writing  $P(X_m \cdots X_M = wuz |  Y_i, Y_{i+1}) = P(\E_u | X_r = u_0, X_s=u_p)  P(X_m \cdots X_r = wu_0, X_s \cdots X_M = u_pz | Y_i, Y_{i+1})$, we obtain
\begin{eqnarray*}
E[ G_i ^2 \, | \, Y_i, Y_{i+1} ] & \geq & \sum_{w \in E^{r}, z \in E^p}  ( G_i(wuz)^2 +  G_i(wvz)^2)  \\
&& \mbox{} \times \epsilon_2 \times
 P(X_m \cdots X_r = wu_0, X_s \cdots X_M = u_pz | Y_i, Y_{i+1}). 
 \end{eqnarray*}
For  $w \in E^{r}$ and $z\in E^p$, observe that
$$|G_i(wuz) - G_i(wvz)| =     \left| \sum_{k=r}^{s-1} (g_k(u_k)) - \sum_{k=r}^{s-1} (g_k(v_k)) \right|  \geq \epsilon,$$
so that
$|G_i(wuz)|^2 + |G_i(wvz)|^2 \geq \epsilon^2$ (since  $a^2+b^2 \geq a^2 +b^2 - 2 |a||b| = (|a| - |b|)^2$). Hence
\begin{eqnarray*}
E[ G_i ^2 \, | \, Y_i, Y_{i+1} ]  &\geq&  \sum_{w \in E^{r}, z \in E^p}  (\epsilon^2 /4) \epsilon_2 P(X_m \cdots X_r = wu_0, X_s \cdots X_M = u_pz | Y_i, Y_{i+1})\\
&  \geq  & (\epsilon^2 /4) \epsilon_2 P( X_r = u_0, X_s= u_p | Y_i, Y_{i+1}) \\
&\geq& \epsilon^2 \epsilon_2\epsilon_1 /4, 
\end{eqnarray*}
and the claim is proved with $\gamma = \epsilon^2 \epsilon_2\epsilon_1 /4$.

\vip
\vip
\vip

{\em Step 2. } To fulfill the proof of Lemma \ref{lem:variance}, we prove that
$$V( \sum_{i=0}^{n-1} f_i(Y_i, \tilde Y_i)) \geq \sum_{i=0}^{n-1} E\left[ \left( f_i(Y_i, \tilde Y_i) - E\left[ f_i(Y_i, \tilde Y_i) | Y_i, Y_{i+1} \right] \right)^2 \right], $$
or, in other words,
$$V( \sum_{i=0}^{n-1} F_i) \geq \sum_{i=0}^{n-1} E\left[ G_i^2 \right]. $$
We recall that $EF_i=0$.
\begin{eqnarray*}
V(\sum_i F_i) &=& E\left[ (\sum_i F_i )^2 \right] \\
&=& E\left[ E\left[ (\sum_i F_i )^2 \, | \, Y_0, \ldots, Y_n \right] \right] \\
&=& E\left[ E\left[ \left(\sum_i (F_i - E[F_i|Y_i,Y_{i+1}]) + \sum_i E[F_i|Y_i,Y_{i+1}] \right)^2 \, | \, Y_0, \ldots, Y_n \right] \right] \\
&=& E\left[ E\left[ \left(\sum_iG_i + \sum_i E[F_i|Y_i,Y_{i+1}] \right)^2 \, | \, Y_0, \ldots, Y_n \right] \right] \\
&=& E\left[ E\left[ \left(\sum_i G_i\right)^2 + H+   \left( \sum_i E[F_i|Y_i,Y_{i+1}] \right)^2  \, | \, Y_0, \ldots, Y_n \right] \right], 
\end{eqnarray*}
where the product term $H$ is
\begin{eqnarray*}
H&=& E\left[    \left(\sum_i G_i\right)  \left( \sum_i E[F_i|Y_i,Y_{i+1}] \right)    \, | \, Y_0, \ldots, Y_n \right] \\
&=&  E\left[  \sum_i G_i   \, | \, Y_0, \ldots, Y_n \right] \left( \sum_i E[F_i|Y_i,Y_{i+1}] \right), 
\end{eqnarray*}
and, for all $i$,
\begin{eqnarray*}
E\left[  G_i   \, | \, Y_0, \ldots, Y_n \right]  &=& E\left[  F_i - E[F_i|Y_i,Y_{i+1}]   \, | \, Y_0, \ldots, Y_n \right] \\
&=& E\left[  F_i \, | \, Y_0, \ldots, Y_n \right]  - E[F_i|Y_i,Y_{i+1}]  \\
&=& 0.
\end{eqnarray*}

To conclude, we observe that, since for all $i <j\leq k $, $E\left[  G_i G_k  \, | \, Y_j \right] = 0$ (by independence of the past and the future conditionally to the present for the Markov chain $(X_n)$),
\begin{eqnarray*}
 E\left[ E\left[ \left(\sum_i G_i \right)^2  \, | \, Y_0, \ldots, Y_n \right] \right] &=&  E\left[ E\left[ \sum_i   \left( G_i \right)^2  \, | \, Y_0, \ldots, Y_n \right] \right]\\
 &=&  \sum_i E\left[   E\left[  \left( G_i)\right)^2  \, | \, Y_0, \ldots, Y_n \right] \right]\\
 &= & \sum_i E\left[   E\left[  \left( G_i)\right)^2  \, | \, Y_i, Y_{i+1} \right] \right] \\
 & \geq & E\left[   \sum_i \gamma \right] \\
 & \geq & n \gamma, 
 \end{eqnarray*}
where the last inequality is indebt to (\ref{var1}).
$\Box$


\section{Specific cases and examples} \label{section:examples}
To get an idea of the whole picture it seems useful  to study what happens in specific cases, making further assumptions on the substitution or/and on the possible choices of $t$. We also give concrete examples at the end of this section.


\subsection{Synchronization}
We consider the case $t=1$ and more generally $t$ integer or say $\dd^{n_0} t$ integer for some $n_0$, i.e.  the expansion of $t$ is finite. To study this simple case, we observe that (at least after a finite number of steps), we can work with a {\it simplified automaton} with state space $\A \times \A$ because when $\tau_k =0$ we never need to know the second letter of $V$ since $m_k+\tau_k \leq d$.

\begin{lem} \label{lem:coboundary} When the expansion of $t$ is finite,
 the subset of states  $\{(a,a) \; ; \; a \in \A\}$ is a recurent class of the (simplified) Markov chain. Moreover, on this class, the cocycle $g$ is a coboundary. Instead of normal fluctuations the limit law of $ \lim_{ n \to \infty} \frac{1}{\sqrt{n}} \int_{0}^{\dd^nt} f \circ h_s(x) ds$ is a  Dirac mass at 0.
 \end{lem}

 \begin{proof}
 We assume that $\tau_{k} = 0$ for all $k$ large enough (larger than $n_{0}$). Thus, the Markov chain is stationary. 
 In the automaton associated to $\tau= 0$, it is clear that  a couple $(a,a)$ is followed by a couple $(b,b)$ since the positions of both  letters are the same in $\sigma(a)$. Therefore the set $\{(a,a) \; ; \; a \in \A\}$ is a recurrent class.

If $\frak{S}^{-n_{0}} x$ belongs  to $a$ and $h_{\dd^{n_0}t}(\frak{S}^{-n_{0}}x)$ is in the same letter $a$ then the integral along the unstable leave does not depend on the exact position of $\frak{S}^{-n_0}x$ in $a$. Indeed shift of the piece of leaf along does not change the value of the integral (what is added from one side is subtracted from the other side).  Moreover, since the value of $\tau_{k} = 0$ is zero, if we shift
$\frak{S}^{-n_{0}} x$ at the beginning of $a$ the point $\frak{S}^{-n_{0}} x$  to $h_{\dd^{n_0}t}(\frak{S}^{-n_{0}}x)$ is also at the beginning of $a$. Thus the leaf from $\frak{S}^{-n_{0}} x$ to $h_{\dd^{n_0}t}(\frak{S}^{-n_{0}}x)$ covers exactly a union of blocks of that scale. Thus, if we refine for all $n > n_0$,  $\int_{0}^{\dd^nt} f(h_s(\sigma^{-n}x)) ds =\int_{0}^{\dd^{n_0}t} f(h_s(\sigma^{-{n_0}}x)) ds$. The value along the Markov chain can change (since it is only an approximation of the integrals, at each scale) but it must keep at bounded distance of the  value of the integral. We claim that the sum of the values along any cycle of the recurrent class $\{(a,a) \; ; \; a \in \A\}$ is zero. In fact,
given a cycle of  class, if the sum of the values of $g$ along this cycle were non zero, we could iterate this cycle a large number of time and get arbitrarily large value for $g$.
Thus $g$ is a coboundary. The limit law of of $ \lim_{ n \to \infty} \frac{1}{\sqrt{n}} \int_{0}^{\dd^nt} f \circ h_s(x) ds$ is therefore a Dirac mass at 0 since  $\int_{0}^{\dd^nt} f \circ h_s(x) ds$ is bounded and we divide by $\sqrt{n}$.
  \end{proof}


\begin{defin}
We say that a letter is {\em synchronizable} if it appears at the same position in the image of two distinct letters. We say that a substitution $\sigma$  is {\em synchronizable} if for every couple of letters $(b,c)$ there is a letter $a$ that appears at the same position in $\sigma(b)$ and $\sigma(c)$. A substitution is  {\em  strongly non synchronizable} if it has no synchronizable letter.
\end{defin}

\begin{prop}
If a substitution $\sigma$ is synchronizable, then the (simplified) Markov chain has a unique recurrent class, the set
$\{(a,a) \; ; \; a \in \A\}$. On this class, the cocycle $g$ is a coboundary. Therefore the limit law of
$ \lim_{ n \to \infty} \frac{1}{\sqrt{n}} \int_{0}^{\dd^nt} f \circ h_s(x) ds$ is a Dirac mass at 0 when the expansion of $t$ in base $\dd$ is finite.
\end{prop}
\begin{proof}
By Lemma \ref{lem:coboundary}, it is enough to prove that the simplified automaton has only one recurrent class.
At each step of the Markov chain there is a positive probability to reach the recurrent class $\{(a,a) \; ; \; a \in \A\}$. Therefore the only recurrent class is the class $\{(a,a) \; ; \; a \in \A\}$.
\end{proof}

\begin{lem}
If $\sub$ has no synchronizable letters, then the recurent class  $\{(a,a) \; ; \; a \in \A\}$ of the (simplified) Markov chain is disconnected of the remainder of the graph.   Hence it is reached  only if reached at initial time.
\end{lem}
\begin{proof} If $b$ and $c$ are different letters of the alphabet $\A$, there is no arrow from the $(b,c)$ to $\{(a,a) \; ; \; a \in \A\}$. Otherwise there would be a synchronizable letter. Therefore the class $\{(a,a) \; ; \; a \in \A\}$ is disconnected from the remainder of the graph.
\end{proof}

In the last paragraph, Figure~\ref{fig:graph1} gives an example of a synchronizable substitution on a 3 letters alphabet.  Figure~\ref{fig:graph2} gives an example of a substitution on a 2 letters alphabet which is non synchronizable.

\subsection{Strongly non synchronizable substitutions on a 2 letters alphabet}
 The following propositions give a recipe to construct strongly non synchronizable substitutions on a 2 letters alphabet.  We also study the ergodic sums for these substitutions for some values of $t$.

\begin{prop} If $\A = \{a,b\}$ a non synchronizable substitution $\sub$ of length $d$ on
$\A$ with eigenvalue 1 is uniquely determined by the image of $a$ ($\sub(b)$ is obtained from $\sub(a)$ by exchanging $a$ by $b$ and vice versa). Moreover the matrix $M$ of $\sub$ has the form
$M = \begin{pmatrix}
k+1 & k \\
k & k+1
\end{pmatrix}
$ where $k +1$ is the number of $a$ in $\sub(a)$ and $d = 2k+1$.
\end{prop}

\begin{proof}
First of all, if $\sub$ is strongly non synchronizable, for every $j$ in $\{1, \dots, d\}$, the letter number $j$ in $\sub(a)$ is different from the letter number $j$ in $\sub(b)$. Therefore one of these letter is an $a$, the other one is a $b$. Thus $\sub(b)$
obtained from $\sub(a)$ is equal to $s (\sub(a))$ where $s$ is the substitution $s(a) = b$, $s(b) = a$. Let $\ell$ be the number of $a$ in $\sub(a)$ and $k$ the number of $b$ in $\sub(a)$. We have $k+ \ell = d$. By the previous discussion, the matrix $M$ has the form
$M = \begin{pmatrix}
\ell & k \\
k & \ell
\end{pmatrix}$. Since
the matrix $M$ has an eigenvalue 1, $\ell = k+1$.
\end{proof}

\begin{prop}
Let $\sub$ be a substitution of length $d$ on  a 2 letters alphabet $\A = \{a,b\}$. We assume that $\sub$ is strongly non synchronizable and that every word of length 3 appears in $\sub(a)$. Then the graph associated to the automaton $\tau =1$ is strongly connected and aperiodic. Moreover $g$ is not a coboundary on this graph. Thus the limit law for $ \lim_{ n \to \infty} \frac{1}{\sqrt{n}} \int_{0}^{\dd^nt} f \circ h_s(x) ds$ is a normal law when the digits of $t$ are ultimately constant equal to 1.
\end{prop}

\begin{rem}
\label{remark-family}
This gives a family of examples that satisfy the hypothesis of Theorem \ref{fluctuations-generic} and hence a proof of Theorem \ref{fluctuations-generic-basic}.
\end{rem}

\begin{proof}
To prove connectedness of the graph,
we first prove that there is an arrow from $(a,aa)$ to any other vertex. Let $(\alpha, \beta\gamma)$ be a triple of letters in $\A$. The word $\alpha \beta\gamma$ appears in $\sub(a)$. We consider the word  $\sub(aa)$. In this word we see
$\alpha W \alpha \beta \gamma$ where the length of $W$ is equal to $d$. Thus the length between $\alpha$ and
$\beta \gamma$ is $d+1$. By definition of the automaton, this means that we have an arrow from $(a,aa)$ to $(\alpha, \beta\gamma)$.

 Now, let us prove that from any triple  $(\alpha, \beta\gamma)$ there is an arrow to $(a,aa)$. If $\alpha = a$ and $\beta =  a$ the word $aaa$ appears somewhere in $\sub(a)$. Thus,
by definition of the automaton with $\tau = 1$, there is an arrow from $(a,a\gamma)$ to $(a,aa)$.
If $\alpha = a$ and $\beta = b$, the word $abb$ appears in $\sub(a)$. At the same position $baa$ appears in $\sub(b)$ since $\sub$ is not strongly synchronizable. Thus there is an arrow from $(a,b\gamma)$ to $(a,aa)$. If $\alpha = b$, $\beta = a$, the same reasonning holds. If $\alpha = \beta = b$, then $aaa$ appears in $\sub(b)$. Thus, by the same argument, there is an arrow from $(b,b\gamma)$ to $(a,aa)$. This proves that the graph is strongly connected.

In fact, we also proved that it is aperiodic since there is an arrow from $(a,a)$ to itself (thus a cycle of length 1).

We prove now that $g$ is not a coboundary. We show that the value of $g$ on some arrow from $(a,aa)$ to itself is different from 0 which is enough. Let $aaa$ be a subword of $\sub(a)$. We have $\sub(aa) = P aaa SP aaa S$ where the length of $SP$ is $d-3$. We consider the arrow from $(a,aa)$ to itself corresponding to $m = \vert P \vert +1$.
The label of this arrow is the number $v_{aaSPa}$. The word $aaSPa$ is a reordering of   $PaaaS$. Thus, the numbers
$v_{aaSPa}$ and $v_{PaaaS}$ coincide. But $v_{PaaaS}  = v_{\sub(a)} = v_{a}$.  One can check that an eigenvector associated to the eigenvalue 1 of $M$ is
$\begin{pmatrix}
1 \\
-1
\end{pmatrix}.$ Thus $v_{a}$ is different from  0.

We now apply Theorem \ref{fluctuations-periodic}. As there is only one recurrent class which is aperiodic, the limit distribution  is a normal law.
\end{proof}

\subsection{Examples} \label{subsection:examples}

In this paragraph, we construct explicit examples with various properties:

\begin{figure}[h]
 \includegraphics[width=8cm]{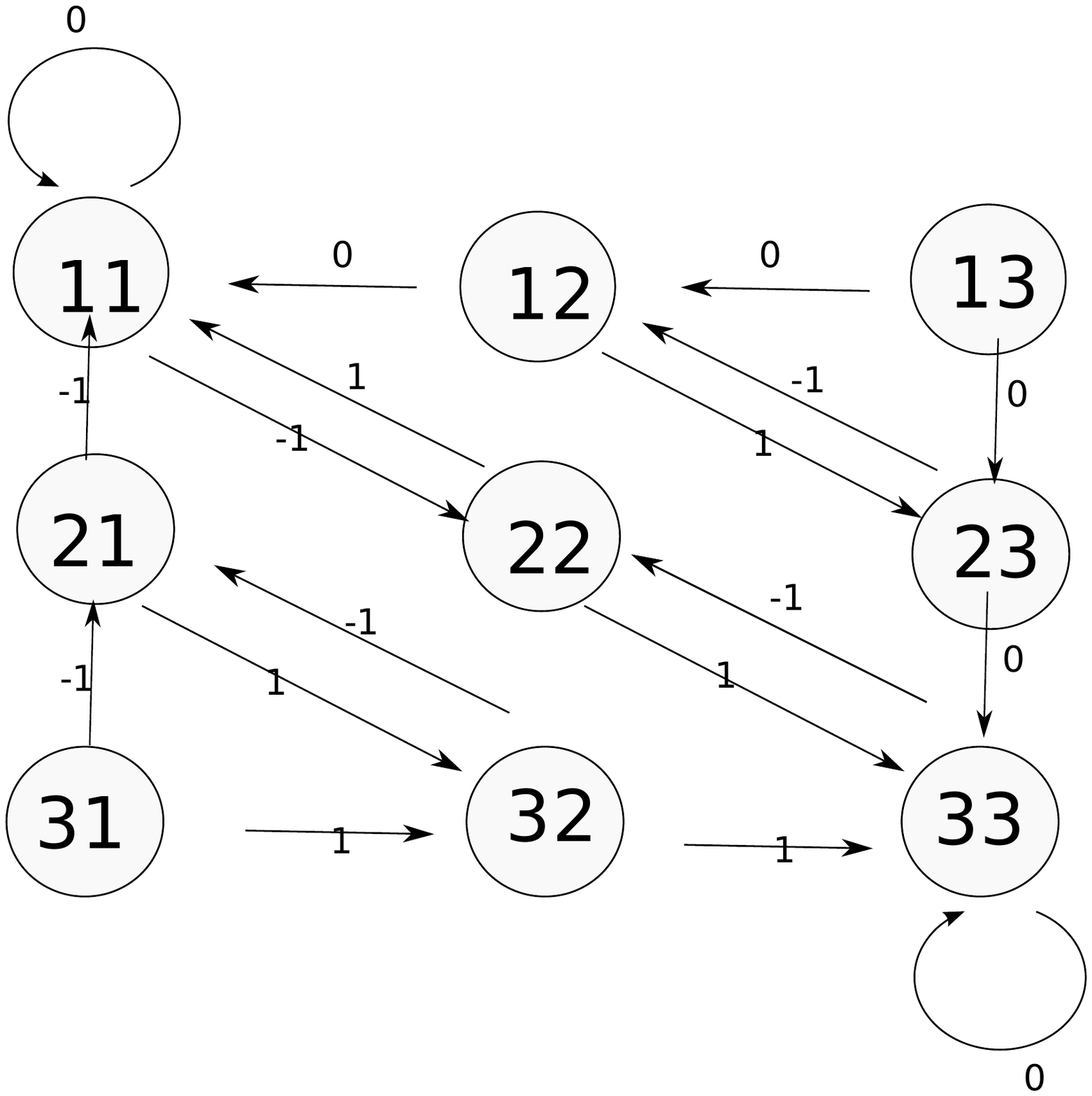}
\caption{\label{fig:graph1}
}
\end{figure}

Figure~\ref{fig:graph1} represents the automatom for
 $$\sigma : \left\{\begin{matrix}
 1 & \to &12 \\
 2 & \to &13\\
 3 & \to & 23
 \end{matrix}\right.$$ and $\tau = 0$.
 The synchronization property is satisfied. The graph has one recurrent component but is not strongly connected (there are transient components). The eigenfunction $f$ is not a coboundary nevertheless, when the 2-adic expansion of $t$ is finite, the cocycle $g$ is a coboundary.

\begin{figure}[h]
 \includegraphics[width=8cm]{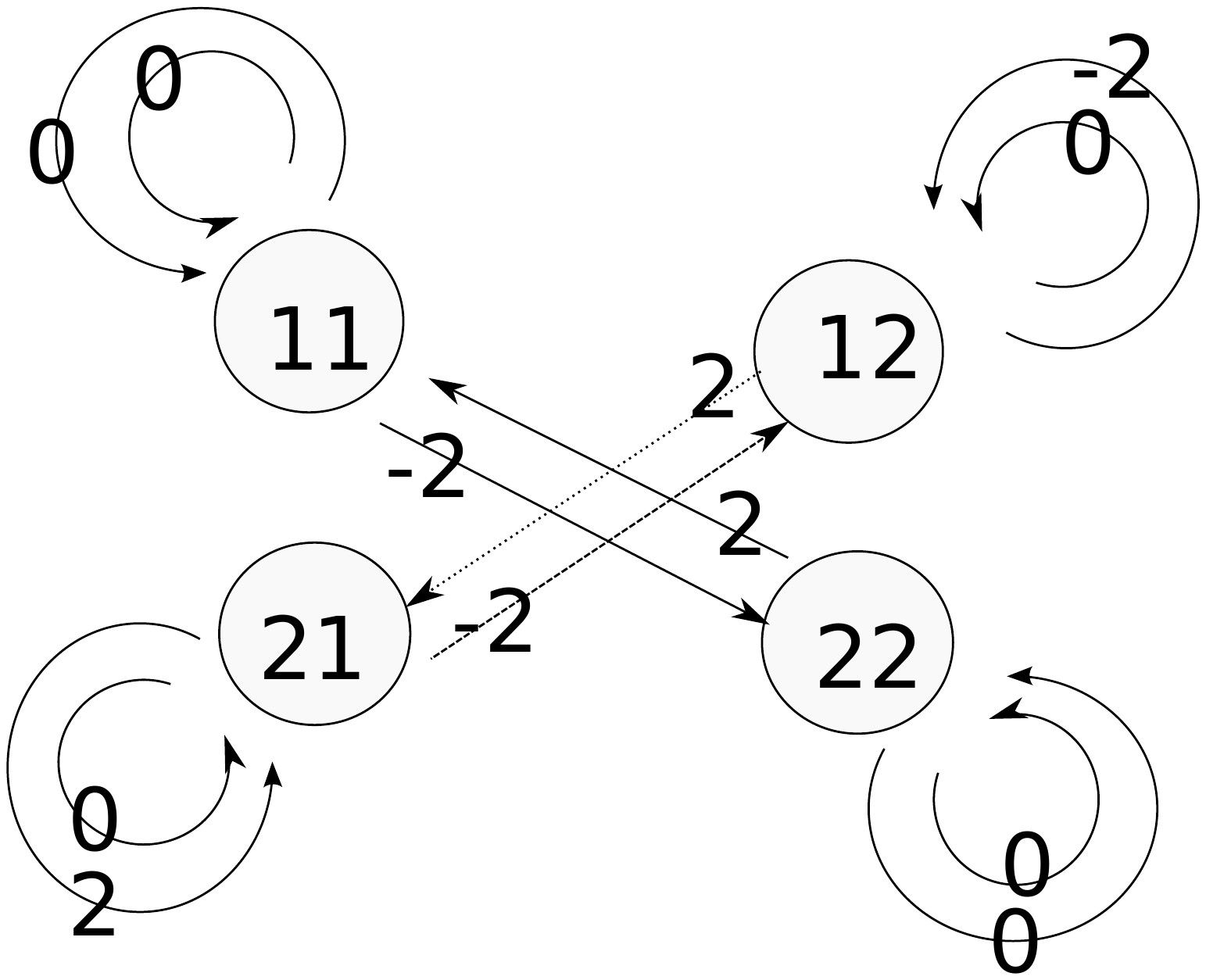}
\caption{\label{fig:graph2}
}
\end{figure}

Figure~\ref{fig:graph2} represents the automaton for
 $$\sigma : \left\{\begin{matrix}
 1 & \to & 112 \\
 2 & \to & 221
 \end{matrix}\right.$$ and $\tau = 0$. The graph has 2 recurrent components:
 the component $\{(a,a) \; ; \; a \in \A\}$ where the cocycle is a coboundary, another component where the Markov chain has non zero variance. Thus, the limiting distribution is a superposition of a Dirac at 0 and a normal distribution. This means that starting from a set of positive measure there is no fluctuation at the scale $\frac{1}{\sqrt{n}}$, on another set of positive measure there are fluctuations at the same scale.

\begin{figure}[p]
 \includegraphics[width=8cm]{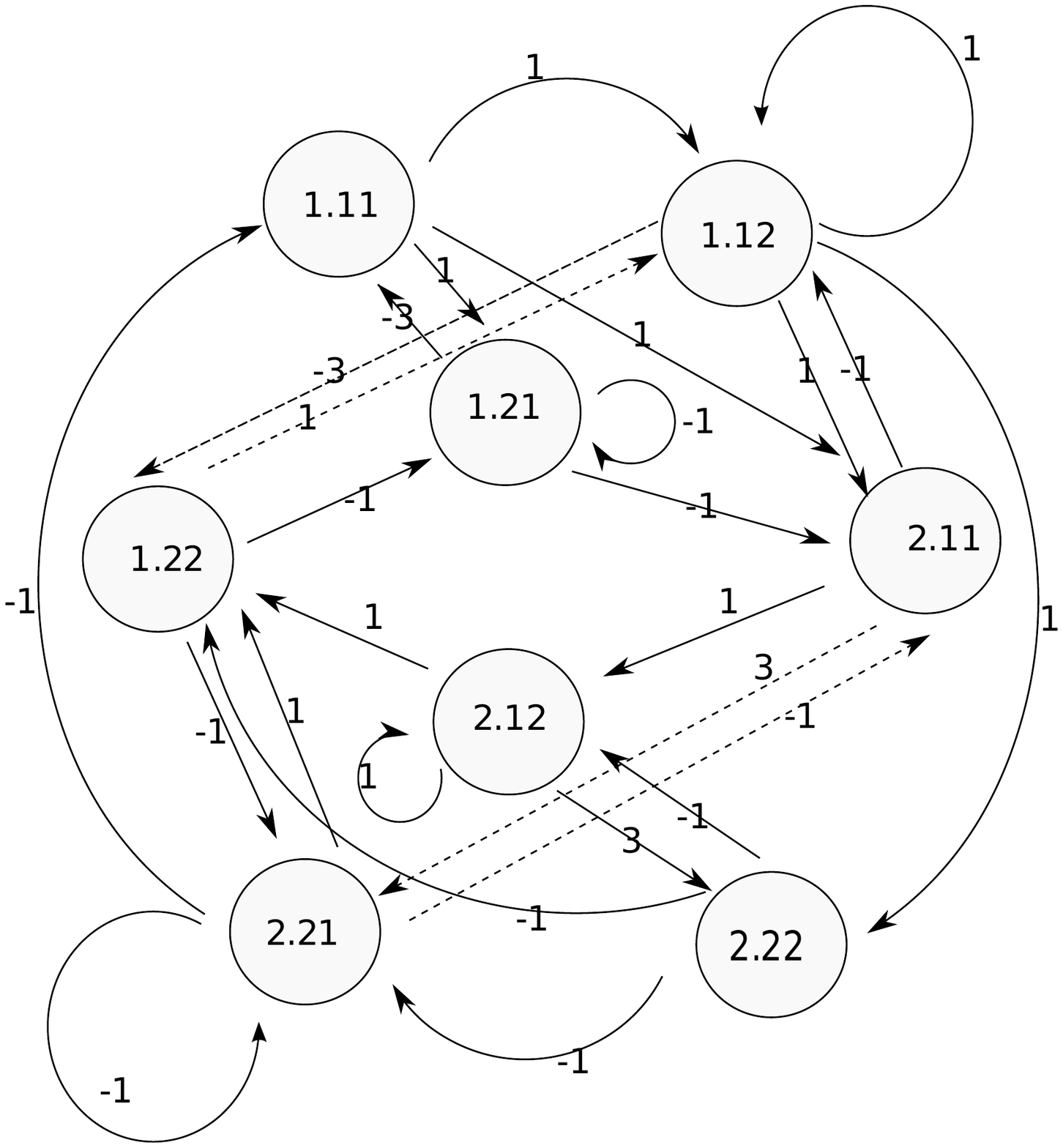}
\caption{\label{fig:graph3}
}
\end{figure}

 Figure~\ref{fig:graph3}  represents the automaton for
 $$\sigma : \left\{\begin{matrix}
 1 & \to & 112 \\
 2 & \to & 221
 \end{matrix}\right.$$ and $\tau = 1$.
 The graph is strongly connected, aperiodic, the variance is positive thus the limiting distribution is a normal distribution.

\begin{figure}[p]
 \includegraphics[width=8cm]{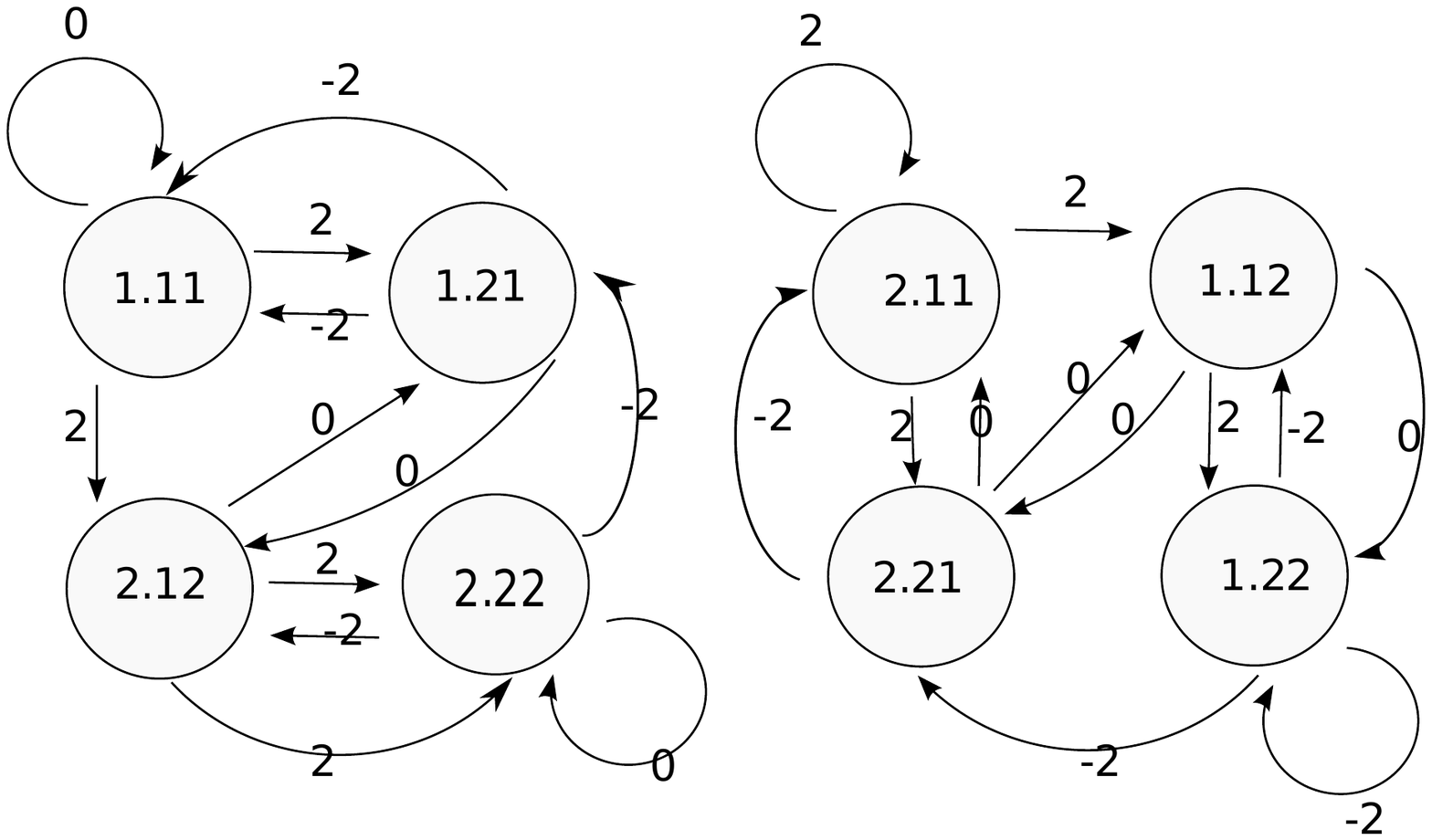}
\caption{\label{fig:graph4}
}
\end{figure}

Figure~\ref{fig:graph4}  represents the automaton for
 $$\sigma : \left\{\begin{matrix}
 1 & \to & 112 \\
 2 & \to & 221
 \end{matrix}\right.
 $$ and $\tau = 2$.

The graph contains two connected components. On one component the cocycle is a coboundary. We have no explanation for this phenomenon. On the other one, there is a positive variance.

\end{document}